\documentclass[10pt]{article}
\usepackage{amstext}
\usepackage[fleqn]{amsmath}
\usepackage{amsfonts}
\usepackage{amssymb}
\usepackage{amsthm}
\usepackage{newlfont}
\usepackage{graphicx}
\usepackage{sectsty}

\theoremstyle{plain}
\theoremstyle{theorem}
\newtheorem{theorem}{Theorem}[section]
\theoremstyle{example}

\theoremstyle{corollary}
\newtheorem{corollary}{Corollary}[section]
\theoremstyle{lemma}

\theoremstyle{proposition}

\theoremstyle{definition}
\newtheorem{definition}{Definition}[section]

\begin{document}

\title{\Large\bf On fractional integrals and derivatives of Bateman's matrix polynomials}

\date{}
\maketitle{}
\author
{\bf Ghazi S. Khammash$^1$,
Shimaa I. Moustafa$^2$,
Shahid Mubeen$^3$
Saralees Nadarajah$^{4, *}$
and Ayman Shehata$^5$}
\\
\\
{\small{\it
$^{1}$Department of Mathematics, Al-Aqsa University, Gaza Strip, Palestine, {\sf ghazikhamash@yahoo.com}
\\
$^{2}$Department of Mathematics, Faculty of Science, Assiut University, Assiut 71516, Egypt, {\sf shimaa1362011@yahoo.com}, {\sf  shimaa\_m@science.aun.edu.eg}
\\
$^{3}$Department of Mathematics, University of Sargodha, Sargodha, Pakistan, {\sf  smjhanda@gmail.com}
\\
$^{4}$Department of Mathematics, University of Manchester, Manchester M13 9PL, UK, {\sf  mbbsssn2@manchester.ac.uk}
\\
$^{5}$Department of Mathematics, Faculty of Science, Assiut University, Assiut 71516, Egypt, {\sf aymanshehata@science.aun.edu.eg}, {\sf drshehata2006@yahoo.com}
\\
$^{*}$corresponding author}.}

\date{}
\maketitle{}

\begin{abstract}
The object of this paper is to investigate certain results involving Bateman's matrix polynomials for integral index.
We obtain some properties of hypergeometric matrix functions, integral representation for generalized hypergeometric
matrix functions and recurrence relations.
We introduce some matrix differential equations of order three,
integral transform and fractional integral formulas for the generalized gauss hypergeometric functions by
using the beta and Laplace transforms formula, Mellin transform of integrable function,
Erdelyi Kober type fractional integral operators and fractional integrals (and fractional derivatives) for the Riemann–Liouville,
and  Weyl operators.
Finally, we present some applications on Bateman's matrix polynomials for different types and Young's matrix functions.
\end{abstract}

\noindent
{\bf   AMS Mathematics Subject Classification  (2020):}  33C99.

\noindent
{\bf  Keywords and phrases:}
Generalized gauss hypergeometric function; Laplace transform; Mellin transform

\section{Introduction}

In recent years, matrix generalization of special functions has grown in importance in
the theory of special functions  \cite{ljc, ljc2}.
One motive is that  special matrix functions solve physical problems and another is that
special matrix functions are closely related to orthogonal matrix polynomials.
Several researchers have studied  special matrix functions,
also many authors have developed numerous matrix integral formulas
involving a variety of special matrix functions in recent years.
The  integrals are of interest in pure mathematics,
are also extremely important in many branches of theoretical and applied physics and engineering.
We refer the readers to \cite{edj, ghgc, ljc, ljc2, ljc4, fmm, hms4}.

Bateman's polynomials relate to the orthogonal polynomials family.
Many  researchers have generalized the classical results on Bateman's polynomials;
and there is a large body of literature,  see  \cite{Rsb1, Rsb2, Rsb4, mpsr, edra, nash, shmr, zsg, hms5}.

In the last two decades, the study of orthogonal matrix polynomials
has been an increasing topic of study, with important results in theory and applications
continuing to appear in the literature  \cite{ggc}.
The theory of orthogonal polynomials also extends to polynomial matrices,
Bateman's matrix polynomials.
We  can see for example \cite{rac3, mme2, mmm, Mam, Ma1, Ma2}.

Let $A=\left(a_{r,s}\right)=A_{1}+i A_{2}=\left(a_{r,s_{1}}\right)+i \left(a_{r,s_{2}}\right)$
be an $m \times n$ array of complex numbers,
$A\in C^{m\times n}$, where $A_{1}=\left(a_{r,s_{1}}\right)$ and $A_{2}=\left(a_{r,s_{2}}\right)$
are real matrices of order $m \times n$,
where $A_{1}, A_{2} \in C^{m,n}$, and $i=\sqrt{-1}$.
$A'$ denotes the transpose of $A$,
and $A^{H}$ denotes the conjugate transpose of $A$,
the $\mathbb{N}$ and $\Bbb{C}$  represent the sets of natural numbers and complex numbers, respectively.
Throughout this paper, for a matrix  $A$ in $\Bbb{C}^{N\times N}$,
its spectrum $\sigma(A)$ denotes the set of all eigenvalues of $A$.
The two-norm will be denoted by $||A||_{2}$ and it is defined by (see \cite{ljc, ljc2})
$||A||_{2}=\sup_{x\neq 0}\frac {||Ax||_{2}}{||x||_{2}}$,
where for a vector $x$ in $\Bbb{C}^{N}$, $||x||_2=\left(x^Tx\right)^\frac {1}{2}$ is the Euclidean norm of $x$.

Let us denote the real numbers $M(A)$ and $m(A)$, respectively, as
\begin{equation}
\displaystyle
M(A)=\max\left\{Re(z): z\in \sigma(A)\right\};
\quad
\displaystyle
m(A)=\min\left\{Re(z): z\in \sigma(A)\right\}.
\label{1.1}
\end{equation}
If $f(z)$ and $g(z)$ are holomorphic functions of
the complex variable $z$ defined in an open set $\Omega$ of the complex plane, and $A$, $B$
are matrices in $\Bbb{C}^{N\times N}$ with $\sigma(A)\subset\Omega$ and $\sigma(B) \subset \Omega$,
such that $AB=BA$, then it follows from the matrix functional calculus properties in \cite{ndjt} that
\begin{equation*}
\displaystyle
f(A)g(B)=g(B)f(A).
\end{equation*}
Throughout this paper, a matrix polynomial of degree $n$ in $x$ is defined as
\begin{equation}
\displaystyle
P_{n}(x)=A_{n}x^{n}+A_{n-1}x^{n-1}+\ldots+A_{1}x+A_{0},
\label{1.2}
\end{equation}
where $x$ is a real variable or complex variable,
$A_j$, $0<j<n$ and $A_{n}\neq \mathbf{0}$ are complex matrices in $\Bbb{C}^{N\times N}$.
The reciprocal gamma function, denoted by $\Gamma^{-1}(z)=\frac {1}{\Gamma(z)}$, is an entire
function of the complex variable $z$, and thus $\Gamma^{-1}(A)$ is a well defined reciprocal
matrix for any matrix $A$ in $\Bbb{C}^{N\times N}$.
Furthermore, if $A$ is a matrix, then
\begin{equation}
\displaystyle
A+nI\quad \text{is\  an invertible matrix\  for\  all\  integers}\  n\geq0,
\label{1.3}
\end{equation}
where $I$ is the identity matrix in $\Bbb{C}^{N\times N}$.
From \cite{ljc, ljc2}, it follows that
\begin{equation*}
\displaystyle
(A)_{n}=A(A+I)\ldots(A+(n-1)I)=\Gamma{(A+nI)}\Gamma^{-1}{(A)};
\quad
\displaystyle
n\geq 1;
\quad
\displaystyle
(A)_{0}=I.
\end{equation*}
If $k$ is large enough, then for $k>|B|$, we  refer to the following relation
that existed in J\'{o}dar and Cort\'{e}s \cite{ljc, ljc2} in the form
\begin{equation*}
\displaystyle
\|(B+kI)^{-1}\|\leq \frac {1}{k-\|B\|};
\quad
\displaystyle
k>\|B\|.
\end{equation*}
If $P$ and $Q$ are positive stable matrices in  $\Bbb{C}^{N\times N}$
such that $Re(z)>0$, $Re(w)>0$, for all $z \in \sigma (P)$, $w \in \sigma (Q)$,
then gamma and beta matrix functions are, respectively, defined as \cite{ljc}
\begin{equation}
\label{1.6}
\displaystyle
\Gamma(P)= \displaystyle \int_{0}^{\infty}e^{-t}t^{P-I}dt
\end{equation}
and
\begin{equation}
\beta(P,Q)=\displaystyle \int_{0}^{1}t^{P-I}(1- t)^{Q-I}dt.
\label{1.7}
\end{equation}
In addition, if $P, Q $ and $P + Q $ are positive stable matrices in $\Bbb{C}^{N\times N}$ and $P$ and $Q$ are commutative, then
\begin{equation}
\displaystyle
\beta(P,Q)=\Gamma(P)\Gamma(Q)\Gamma^{-1}(P+Q).
\label{1.8}
\end{equation}
If $A(k,n)$ and $B(k,n)$ are matrices in $\Bbb{C}^{N\times N}$  for $n\geq 0$ and $k\geq 0$,
then it follows, in a manner analogous to the proof of Lemma \textbf{11} \cite{edj}, that
\begin{equation}
\begin{split}
\displaystyle
\sum_{n=0}^{\infty}\sum_{\ell=0}^{\infty}A(k,n)=\sum_{n=0}^{\infty}\sum_{\ell=0}^{\left[\frac {1}{2}n\right]}A(k,n-2k),
\\
\displaystyle
\sum_{n=0}^{\infty}\sum_{\ell=0}^{\infty}B(k,n)=\sum_{n=0}^{\infty}\sum_{\ell=0}^{n}B(k,n-k).
\end{split}
\label{1.10}
\end{equation}
Similarly to (\ref{1.10}), we can write
\begin{equation*}
\begin{split}
\displaystyle
\sum_{n=0}^{\infty}\sum_{\ell=0}^{\left[\frac {1}{2}n\right]}A(l,n)=\sum_{n=0}^{\infty}\sum_{\ell=0}^{\infty}A(l,n+2l),
\\
\displaystyle
\sum_{n=0}^{\infty}\sum_{\ell=0}^{n}B(k,n)=\sum_{n=0}^{\infty}\sum_{\ell=0}^{\infty}B(l,n+l).
\end{split}
\end{equation*}
The hypergeometric matrix function  ${}_{2}F_{1}(A,B;C;z)$ is given in the form
\begin{equation*}
\displaystyle
{}_{2}F_{1}(A,B;C;z)=\sum_{\ell=0}^{\infty}\frac {(A)_l(B)_l\left[(C)_l\right]^{-1}}{l!}z^l
\end{equation*}
for $A$, $B$, and $C$ matrices in $\Bbb{C}^{N\times N}$ such that $C+nI$ is an invertible
matrix for all integers $n\geq 0$ and for $|z|<1$.
J\'{o}dar and Cort\'{e}s \cite{ljc2, ljc4} observed that the series is absolutely
convergent for $|z|=1$ when $m(C)> M(A)+M(B)$, where $m(P)$ and $M(P)$ are  in (\ref{1.1}) for any matrix $P$ in $\Bbb{C}^{N\times N}$.

\begin{definition}
\label{dfn 1.1}
Let $f$ be an entire function of $z$ in the form (see \cite{bo,co,  ho, rli})
\begin{eqnarray*}
\displaystyle
f(z)=\sum_{s=0}^{\infty}a_{s}z^{s}
\end{eqnarray*}
then, we define the order $\rho$ and type $\tau$ of the function $f(z)$ by
\begin{eqnarray*}
\displaystyle
\rho(f)=\limsup_{s \to \infty}\frac {s\ln(s)}{\ln\left(\frac {1}{\left |a_{s}\right  |}\right)}
\end{eqnarray*}
and
\begin{eqnarray*}
\displaystyle
\tau=\frac {1}{e\rho}\limsup_{s \to \infty}s\left(\left |a_{s} \right  |\right)^{\frac {\rho}{s}}.
\end{eqnarray*}
\end{definition}
A study of the order and type of hypergeometric matrix functions ${}_{1}F_{2}$
and Bateman's matrix polynomials serves as a link between the subjects of special functions and complex analysis.
We establish  important properties of hypergeometric matrix functions ${}_{1}F_{2}$ and Bateman's matrix polynomials.

\section{Properties of hypergeometric matrix functions ${}_{1}F_{2}$ and main results}

In this section, we give the definition of  ${}_{1}F_{2}$.

\begin{definition}
Let $A$, $B$ and $C$ be commutative matrices in $\Bbb{C}^{N\times N}$ satisfying
\begin{eqnarray}
\label{2.1}
\displaystyle
B+sI\  \text{and}\  C+sI\quad \text{are\  invertible matrices}\  \forall s\in \Bbb{N}\cup\{0\}.
\end{eqnarray}
Then we define  ${}_{1}F_{2}$ as
\begin{eqnarray}
\displaystyle
{}_{1}F_{2}(A;B,C;z)=\sum_{s=0}^{\infty}\frac {z^{s}}{s!}(A)_{s}\left[(B)_{s}\right]^{-1}\left[(C)_{s}\right]^{-1}=\sum_{s=0}^{\infty}z^{s}U_{s}.
\label{2.2}
\end{eqnarray}
\end{definition}

\begin{theorem}
\label{thm 2.1}
If the condition in (\ref{2.2}) holds, then
\begin{eqnarray}
\label{2.3}
\displaystyle
{}_{1}F_{2}(A;B,C;z)=\sum_{s=0}^{\infty}z^{s}U_{s}
\end{eqnarray}
is an entire function of $z$.
\end{theorem}

\begin{proof}
Now, we find the radius of convergence $R$ of ${}_{1}F_{2}$.
With the help of Cauchy-Hadamard formula in \cite{Kams} and (\ref{2.2}),
\begin{align*}
\displaystyle
\frac {1}{R}
&=
\displaystyle
\limsup_{s\rightarrow\infty}
\left(\|U_{s}\|\right)^{\frac {1}{s}}=
\lim_{s\to\infty}\sup\left(\bigg{\|}\frac {(A)_{s}\left[(B)_{s}\right]^{-1}
\left[(C)_{s}\right]^{-1}}{s!}\bigg{\|}\right)^{\frac {1}{s}}
\\
&=
\displaystyle
\limsup_{s\rightarrow\infty}\bigg{\|}\frac {\Gamma(A+sI)\Gamma^{-1}(A)\Gamma^{-1}(B+sI)\Gamma(B)\Gamma^{-1}(C+sI)\Gamma(C)}{s!}\bigg{\|}^{\frac {1}{s}}
\\
&=
\displaystyle
\limsup_{s\rightarrow\infty}\bigg{\|}\sqrt{2\pi}e^{-(A+sI)}(A+sI)^{A+(s-\frac {1}{2})I}\bigg{(}\sqrt{2\pi}e^{-(B+sI)}
(B+sI)^{B+ \left(s-\frac {1}{2}\right)I}\bigg{)}^{-1}
\\
&
\displaystyle
\bigg{(}\sqrt{2\pi}e^{-(C+sI)}(C+sI)^{C+(s-\frac {1}{2})I}\bigg{)}^{-1}
\frac {\Gamma^{-1}(A)\Gamma(B)\Gamma(C)}{\sqrt{2\pi}e^{-s}s^{s+\frac {1}{2}}}\bigg{\|}^{\frac {1}{s}}
\\
&=
\displaystyle
\limsup_{s\rightarrow\infty}\bigg{\|}\sqrt{2\pi}e^{-(A+sI)}(A+sI)^{A+ \left(s-\frac {1}{2}\right)I}
\frac {1}{\sqrt{2\pi}}e^{(B+sI)}(B+sI)^{-B-\left(s-\frac {1}{2}\right)I}
\\
&
\quad
\displaystyle
\frac {1}{\sqrt{2\pi}}e^{(C+sI)}(C+sI)^{-C-\left(s+\frac {1}{2}\right)I}
\frac {\Gamma^{-1}(A)\Gamma(B)\Gamma(C)}{\sqrt{2\pi}e^{-s}s^{s+\frac {1}{2}}}\bigg{\|}^{\frac {1}{s}}
\\
&
\displaystyle
\thickapprox
\limsup_{s\rightarrow\infty}\bigg{\|}e^{-(A+sI)}(A+sI)^{A+\left(s-\frac {1}{2}\right)I}e^{(B+sI)}(B+sI)^{-B-\left(s-\frac {1}{2}\right)I}
\\
&
\quad
\displaystyle
e^{(C+sI)}(C+sI)^{-C-(s-\frac {1}{2})I}\frac {1}{e^{-s}s^{s+\frac {1}{2}}}\bigg{\|}^{\frac {1}{s}}
\\
&
\displaystyle
\thickapprox
\limsup_{s\rightarrow\infty}\bigg{\|}e^{B+sI+C+sI-A-sI+sI}(A+sI)^{A+\left(s-\frac {1}{2}\right)I}(B+sI)^{-B-\left(s-\frac {1}{2}\right)I}
\\
&
\quad
\displaystyle
(C+sI)^{-C-(s-\frac {1}{2})I}s^{-s-\frac {1}{2}}\bigg{\|}^{\frac {1}{s}}
\\
&
\displaystyle
\thickapprox
\limsup_{s\rightarrow\infty}\bigg{\|}e^{B+C-A+2sI}(A+sI)^{A+\left(s-\frac {1}{2}\right)I}(B+sI)^{-B-\left(s-\frac {1}{2}\right)I}
\\
&
\quad
\displaystyle
(C+sI)^{-C-\left(s-\frac {1}{2}\right)I}s^{-s-\frac {1}{2}}\bigg{\|}^{\frac {1}{s}}
\\
&
\displaystyle
\thickapprox
\limsup_{s\rightarrow\infty}
\bigg{\|}e^{B+C-A}e^{2s}(A+sI)^{A+\left(s-\frac {1}{2}\right)I}
(B+sI)^{-B-\left(s-\frac {1}{2}  \right)I}
(C+sI)^{-C-\left(s-\frac {1}{2}\right)I}
s^{-s-\frac {1}{2}}\bigg{\|}^{\frac {1}{s}}
\\
&
\displaystyle
\thickapprox
e^{2}\limsup_{s\rightarrow\infty}
\bigg{\|}(A+sI)^{A+\left(s-\frac {1}{2}\right)I}
(B+sI)^{-B-\left(s-\frac {1}{2}\right)I}(C+sI)^{-C-\left(s-\frac {1}{2}\right)I}s^{-s-\frac {1}{2}}\bigg{\|}^{\frac {1}{s}}
\\
&
\displaystyle
\thickapprox
e^{2}\limsup_{s\rightarrow\infty}
\bigg{\|}(A+sI)^{A-\frac {1}{2}I}(A+sI)^{s}(B+sI)^{-B+\frac {1}{2}I}
(B+sI)^{-s}(C+sI)^{-C+\frac {1}{2}I}(C+sI)^{-s}s^{-s-\frac {1}{2}}\bigg{\|}^{\frac {1}{s}}
\\
&
\displaystyle
\thickapprox
e^{2}\limsup_{s\rightarrow\infty}
\bigg{\|}
\frac {(A+sI)(B+sI)^{-1}(C+sI)^{-1}}{s}\bigg{\|}\bigg{\|}(A+sI)^{A-\frac {1}{2}I}(B+sI)^{-B+\frac {1}{2}I}
(C+sI)^{-C+\frac {1}{2}I}s^{-\frac {1}{2}}\bigg{\|}^{\frac {1}{s}}=0.
\end{align*}
Thus, (\ref{2.2}) is convergent for all complex numbers $z$.
That is, ${}_{1}F_{2}$ is an entire function under  (\ref{2.1}) of the definition (\ref{2.2}).
\end{proof}

\begin{theorem}
\label{thm 2.2}
If the conditions stated in Theorem 2.1 hold, then
\begin{align*}
\displaystyle
{}_{1}F_{2}=\sum_{s=0}^{\infty}z^{s}U_{s}
\end{align*}
is an entire function of order $\frac {1}{2}$.
\end{theorem}

\begin{proof}
We can  calculate the order of ${}_{1}F_{2}$ of the complex variable as follows
\begin{eqnarray}
\label{2.5}
\begin{split}
\displaystyle
\rho\left({}_{1}F_{2}\right)
&=
\displaystyle
\limsup_{s \to
\infty}\bigg{\|}\frac {s\ln(s)}{\ln\left(\frac {1}{U_{s}}\right)}\bigg{\|}=
\limsup_{s \to \infty}\bigg{\|}\frac {s\ln(s)}{\ln\left(s!\ (B)_{s}(C)_{s}\left[(A)_{s}\right]^{-1}\right)}\bigg{\|}
\\
&=
\displaystyle
\limsup_{s \to \infty}\bigg{\|}\frac {s\ln(s)}{\ln\left(s!\ \Gamma(B+sI)\Gamma^{-1}(B)\Gamma(C+sI)\Gamma^{-1}(C)\Gamma^{-1}(A+sI)\Gamma(A)\right)}\bigg{\|}
\\
&=
\displaystyle
\limsup_{k \to \infty}\bigg{\|}\frac {1}{\Psi}\bigg{\|}=\limsup_{k \to \infty}\bigg{\|}
\frac {1}{0+0+I-0+0+I-0+0+I-0+0-I-0}\bigg{\|}=\frac {1}{2},
\end{split}
\end{eqnarray}
where
\begin{align*}
\displaystyle
\Psi
&=
\displaystyle
\frac {\ln\Gamma(A)-\ln\Gamma(B)-\ln\Gamma(C)}{s\ln(s)}+\frac {1}{2}\frac {\ln(2\pi\ s)}{s\ln(s)}I+\frac {s\ln(s)}{s\ln(s)}I-\frac {s\ln(e)}{s\ln(s)}I
\\
&
\quad
\displaystyle
+\frac {1}{2}\frac {\ln\left(2\pi(B+(s-1)I)\right)}{s\ln(s)}+
\frac {(B+(s-1)I)\ln\left(B+(s-1)I\right)}{s\ln(s)}-\frac {\left(B+(s-1)I\right)\ln(e)}{s\ln(s)}
\\
&
\quad
\displaystyle
+\frac {1}{2}\frac {\ln\left(2\pi(C+(s-1)I)\right)}{s\ln(s)}+
\frac {\left(C+(s-1)I\right)\ln\left(C+(s-1)I\right)}{s\ln(s)}-\frac {\left(C+(s-1)I\right)\ln(e)}{s\ln(s)}
\\
&
\quad
\displaystyle
-\frac {1}{2}\frac {\ln\left(2\pi \left(A+(s-1)I\right)\right)}{s\ln(s)}-
\frac {\left(A+(s-1)I\right)\ln\left(A+(s-1)I\right)}{s\ln(s)}+\frac {\left(A+(s-1)I\right)\ln(e)}{s\ln(s)}.
\end{align*}
The proof is complete.
\end{proof}

\begin{theorem}
\label{thm 2.3}
If the condition (\ref{2.2}) holds, then the type of function
\begin{eqnarray}
\label{2.6}
\displaystyle
{}_{1}F_{2}(A;B,C;z)=\sum_{s=0}^{\infty}z^{s}U_{s}
\end{eqnarray}
is an entire function of $2$.
\end{theorem}

\begin{proof}
The hypergeometric matrix function ${}_{1}F_{2}$ can be formulated as
\begin{align*}
\displaystyle
\tau
&=
\displaystyle
\tau\left( {}_{1}F_{2}\right)=\frac {1}{e\rho}\limsup_{s \to \infty}
\bigg{\|}s\left(U_{s}\right)^{\frac {\rho}{s}}\bigg{\|}=
\frac {1}{e\rho}\limsup_{s \to \infty}\bigg{\|}s\bigg(\frac {(A)_{s}
\left[(B)_{s}\right]^{-1}\left[(C)_{s}\right]^{-1}}{s!}\bigg)^{\frac {\rho}{s}}\bigg{\|}
\\
&=
\displaystyle
\frac {1}{e\rho}\limsup_{s\rightarrow\infty}s\bigg{\|}
\frac {\Gamma(A+sI)\Gamma^{-1}(A)\Gamma^{-1}(B+sI)\Gamma(B)\Gamma^{-1}(C+sI)\Gamma(C)}{s!}\bigg{\|}^{\frac {\rho}{s}}
\\
&=
\displaystyle
\frac {1}{e\rho}\limsup_{s\rightarrow\infty}s\bigg{\|}\sqrt{2\pi}e^{-(A+sI)}(A+sI)^{A+sI-\frac {1}{2}I}
\bigg{(}\sqrt{2\pi}e^{-(B+sI)}(B+sI)^{B+sI-\frac {1}{2}I}\bigg{)}^{-1}
\\
&
\quad
\displaystyle
\bigg{(}\sqrt{2\pi}e^{-(C+sI)}(C+sI)^{C+sI-\frac {1}{2}I}\bigg{)}^{-1}
\frac {\Gamma^{-1}(A)\Gamma(B)\Gamma(C)}{\sqrt{2\pi}e^{-s}s^{s+\frac {1}{2}}}\bigg{\|}^{\frac {\rho}{s}}
\\
&=
\displaystyle
\frac {1}{e\rho}\limsup_{s\rightarrow\infty}s\bigg{\|}\sqrt{2\pi}e^{-(A+sI)}(A+sI)^{A+sI-\frac {1}{2}I}
\frac {1}{\sqrt{2\pi}}e^{(B+sI)}(B+sI)^{-B-sI+\frac {1}{2}I}
\\
&
\quad
\displaystyle
\frac {1}{\sqrt{2\pi}}e^{(C+sI)}(C+sI)^{-C-sI+\frac {1}{2}I}
\frac {\Gamma^{-1}(A)\Gamma(B)\Gamma(C)}{\sqrt{2\pi}e^{-s}s^{s+\frac {1}{2}}}\bigg{\|}^{\frac {\rho}{s}}
\\
&
\displaystyle
\thickapprox
\frac {1}{e\rho}\limsup_{s\rightarrow\infty}s\bigg{\|}e^{-(A+sI)}(A+sI)^{A+sI-\frac {1}{2}I}e^{(B+sI)}(B+sI)^{-B-sI+\frac {1}{2}I}
\\
&
\quad
\displaystyle
e^{(C+sI)}(C+sI)^{-C-sI+\frac {1}{2}I}
\frac {1}{e^{-s}s^{s+\frac {1}{2}}}\bigg{\|}^{\frac {\rho}{s}}
\\
&
\displaystyle
\thickapprox
\frac {1}{e\rho}\limsup_{s\rightarrow\infty}s\bigg{\|}e^{B+sI+C+sI-A-sI+sI}(A+sI)^{A+sI-\frac {1}{2}I}
(B+sI)^{-B-sI+\frac {1}{2}I}
\\
&
\quad
\displaystyle
(C+sI)^{-C-sI+\frac {1}{2}I}s^{-s-\frac {1}{2}}\bigg{\|}^{\frac {\rho}{s}}
\\
&
\displaystyle
\thickapprox
\frac {1}{e\rho}\limsup_{s\rightarrow\infty}s\bigg{\|}e^{B+C-A+2sI}(A+sI)^{A+sI-\frac {1}{2}I}
(B+sI)^{-B-sI+\frac {1}{2}I}
\\
&
\displaystyle
\quad
(C+sI)^{-C-sI+\frac {1}{2}I}s^{-s-\frac {1}{2}}\bigg{\|}^{\frac {\rho}{s}}
\\
&
\displaystyle
\thickapprox\frac {1}{e\rho}\limsup_{s\rightarrow\infty}s\bigg{\|}
e^{B+C-A}e^{2s}(A+sI)^{A+sI-\frac {1}{2}I}(B+sI)^{-B-sI+\frac {1}{2}I}(C+sI)^{-C-sI+\frac {1}{2}I}
s^{-s-\frac {1}{2}}\bigg{\|}^{\frac {\rho}{s}}
\\
&
\displaystyle
\thickapprox \frac {1}{e\rho}
e^{2\rho}\limsup_{s\rightarrow\infty}s\bigg{\|}(A+sI)^{A+sI-\frac {1}{2}I}
(B+sI)^{-B-sI+\frac {1}{2}I}(C+sI)^{-C-sI+\frac {1}{2}I}s^{-s-\frac {1}{2}}\bigg{\|}^{\frac {\rho}{s}}
\\
&
\displaystyle
\thickapprox
\frac {1}{e\rho} e^{2\rho}\limsup_{s\rightarrow\infty}s\bigg{\|}(A+sI)^{A-\frac {1}{2}I}
(A+sI)^{s}(B+sI)^{-B+\frac {1}{2}I}(B+sI)^{-s}(C+sI)^{-C+\frac {1}{2}I}(C+sI)^{-s}s^{-s-\frac {1}{2}}\bigg{\|}^{\frac {\rho}{s}}
\\
&
\displaystyle
\thickapprox
\frac {1}{e\rho} e^{2\rho}
\limsup_{s\rightarrow\infty}s\bigg{\|}\frac {(A+sI)(B+sI)^{-1}(C+sI)^{-1}}{s}\bigg{\|}^{\rho}\bigg{\|}
(A+sI)^{A-\frac {1}{2}I}(B+sI)^{-B+\frac {1}{2}I}(C+sI)^{-C+\frac {1}{2}I}s^{-\frac {1}{2}}\bigg{\|}^{\frac {\rho}{s}}
\\
&
\displaystyle
\thickapprox
\frac {1}{e\rho} e^{2\rho}=\frac {e^{2\rho-1}}{\rho}=\frac {1}{\rho}=2.
\end{align*}
The proof is complete.
\end{proof}

Next, we derive an integral representation for ${}_{1}F_{2}$ by proving the following result.

\begin{theorem}
\label{thm 2.4}
Let $A$, $B$ and $C$ be matrices in $\Bbb{C}^{N\times N}$ and $Re(A)>0$, $Re(B)>0$, $Re(C)>0$, $A_{i}$,
$Re\left(A_{i}\right)>0$, $1\leq i\leq r$ and $B_{j}$, $Re\left(B_{j}\right)>0$, $1\leq j\leq s$
such that $B-A$ and $C-A$ are staple positive invertible matrices,
where $r$ and $s$ are finite positive integers.
Then we have the integral representation
\begin{align}
\displaystyle
{}_{1}F_{2}(A;B,C;z)
&=
\displaystyle
\Gamma(B)\Gamma^{-1}(A)\Gamma^{-1}(B-A)\int_{0}^{1}t^{A-1}(1-t)^{B-A-1}\ {}_{0}F_{1}(-;C;zt)dt
\nonumber
\\
&=
\displaystyle
\Gamma(C)\Gamma^{-1}(A)\Gamma^{-1}(C-A)\int_{0}^{1}t^{A-1}(1-t)^{C-A-1}\ {}_{0}F_{1}(-;B;zt)dt
\label{2.7}
\end{align}
for all $z$ in the complex plane cut along the real axis from $1$ to $\infty$.
Here, it is understood that $\arg(t)=arg(1-t)=0$ and $(1-zt)^{-a}$ has its principal value.
\end{theorem}

\begin{proof}
From  (\ref{1.7}) and (\ref{1.8}), we have
\begin{align*}
\displaystyle
\int_{0}^{1}t^{A+(k-1)I}(1-t)^{B-A-1}dt
&=
\displaystyle
\mathbf{B}(A+kI,B-A)=\Gamma(A+kI)\Gamma(B-A)\Gamma^{-1}(B+kI)
\\
&=
\displaystyle
\Gamma(A)\Gamma^{-1}(A)\Gamma(A+kI)\Gamma(B-A)\Gamma^{-1}(B)\Gamma(B)\Gamma^{-1}(B+kI)
\\
&=
\displaystyle
\Gamma(A)\Gamma(B-A)\Gamma^{-1}(B)(A)_{k}\left[(B)_{k}\right]^{-1}.
\end{align*}
Similarly, we can establish (\ref{2.7}) by using (\ref{1.7}) and (\ref{1.8}) as follows
\begin{align*}
\displaystyle
\int_{0}^{1}t^{A+(k-1)I}(1-t)^{C-A-1}dt
&=
\displaystyle
\mathbf{B}(A+kI,C-A)=\Gamma(A+kI)\Gamma(C-A)\Gamma^{-1}(C+kI)
\\
&=
\displaystyle
\Gamma(A)\Gamma^{-1}(A)\Gamma(A+kI)\Gamma(C-A)\Gamma^{-1}(C)\Gamma(C)\Gamma^{-1}(C+kI)
\\
&=
\displaystyle
\Gamma(A)\Gamma(C-A)\Gamma^{-1}(C)(A)_{k}\left[(C)_{k}\right]^{-1}.
\end{align*}
The proof is complete.
\end{proof}

\begin{theorem}
\label{thm 2.5}
Let  $A$, $B$, $A_{1}, B_{1}$ and $C$ be commuting positive stable
matrices in $\Bbb{C}^{N\times N}$ and  $B_{1}-A_{1}$ a staple positive invertible matrix.
Then, for $|z|<1$, we have
\begin{eqnarray}
\displaystyle
\int_{0}^{1}t^{A_{1}-1}(1-t)^{B_{1}-A_{1}-1} \ {}_{1}F_{2}(A;B,C;zt)dt=
\Gamma\left(A_{1}\right)
\Gamma\left(B_{1}-A_{1}\right)\Gamma^{-1}\left(B_{1}\right)\ {}_{2}F_{3}\left(A,A_{1};B,B_{1},C;z\right)
\label{2.9}
\end{eqnarray}
and
\begin{eqnarray}
\displaystyle
\int_{0}^{1}t^{A_{1}-1}(1-t)^{B_{1}-1}\  {}_{1}F_{2}(A;B,C;zt)dt=
\Gamma\left(A_{1}\right)\Gamma\left(B_{1}\right)
\Gamma^{-1}\left(A_{1}+B_{1}\right)
\ {}_{2}F_{3}\left(A,A_{1};B,B+B_{1},C;z\right).
\label{2.10}
\end{eqnarray}
\end{theorem}

\begin{proof}
For (\ref{2.9}), using (\ref{1.7}), (\ref{1.8}) and (\ref{2.2}), we have
\begin{eqnarray*}
\displaystyle
\int_{0}^{1}t^{A_{1}+(k-1)I}(1-t)^{B_{1}-A_{1}-I}dt=
\mathbf{B}\left(A_{1}+kI,B_{1}-A_{1}\right)=\Gamma\left(A_{1}\right)
\Gamma\left(B_{1}-A_{1}\right)\Gamma^{-1}\left(B_{1}\right)
\left(A_{1}\right)_{k}\left[\left(B_{1}\right)_{k}\right]^{-1}
\end{eqnarray*}
and
\begin{align*}
&
\displaystyle
\sum_{k=0}^{\infty}\frac {(A)_{k}}{k!}\left[(B)_{k}\right]^{-1}
\left[(C)_{k}\right]^{-1}z^{k}\int_{0}^{1}t^{A_{1}+(k-1)I}(1-t)^{B_{1}-A_{1}-1}dt
\\
&=
\displaystyle
\Gamma\left(A_{1}\right)
\Gamma\left(B_{1}-A_{1}\right)
\Gamma^{-1}\left(B_{1}\right)
\sum_{k=0}^{\infty}\frac {(A)_{k}\left(A_{1}\right)_{k}}{k!}\left[(B)_{k}\right]^{-1}\left[(C)_{k}\right]^{-1}\left[\left(B_{1}\right)_{k}\right]^{-1}z^{k}
\\
&=
\displaystyle
\Gamma\left(A_{1}\right)
\Gamma\left(B_{1}-A_{1}\right)
\Gamma^{-1}\left(B_{1}\right)\ {}_{2}F_{3}\left(A,A_{1};B,B_{1},C;z\right).
\end{align*}
Now for (\ref{2.10}) using (\ref{1.7}), (\ref{1.8})    and (\ref{2.2}), we obtain
\begin{eqnarray*}
\displaystyle
\int_{0}^{1}t^{A_{1}+(k-1)I}(1-t)^{B_{1}-I}dt=
\mathbf{B}\left(A_{1}+kI,B_{1}\right)=\Gamma\left(A_{1}\right)\Gamma\left(B_{1}\right)
\Gamma^{-1}\left(A_{1}+B_{1}\right)
\left(A_{1}\right)_{k}
\left[\left(A_{1}+B_{1}\right)_{k}\right]^{-1}
\end{eqnarray*}
and
\begin{align*}
&
\displaystyle
\sum_{k=0}^{\infty}\frac {(A)_{k}}{k!}\left[(B)_{k}\right]^{-1}\left[(C)_{k}\right]^{-1}z^{k}\int_{0}^{1}t^{A_{1}+(k-1)I}(1-t)^{B_{1}-1}dt
\\
&=
\displaystyle
\Gamma\left(A_{1}\right)\Gamma\left(B_{1}\right)
\Gamma^{-1}\left(A_{1}+B_{1}\right)
\sum_{k=0}^{\infty}\frac {(A)_{k}\left(A_{1}\right)_{k}}{k!}\left[(B)_{k}\right]^{-1}\left[(C)_{k}\right]^{-1}
\left[\left(A_{1}+B_{1}\right)_{k}\right]^{-1}z^{k}
\\
&=
\displaystyle
\Gamma\left(A_{1}\right)\Gamma\left(B_{1}\right)
\Gamma^{-1}\left(A_{1}+B_{1}\right)\ {}_{2}F_{3}\left(A,A_{1};B,B+B_{1},C;z\right).
\end{align*}
Interpreting the right hand side of (\ref{2.3}), in the view of  (\ref{1.7}) and the
concept of the Hadamard given by (\ref{1.3}), we arrive at the required result (\ref{2.10}).
\end{proof}

\begin{theorem}
\label{thm 2.6}
Let  $A$, $B$, $P$, $Q$ and $C$ be commuting positive stable matrices in $\Bbb{C}^{N\times N}$.
Then, for $|z|<1$, we have the following representations
\begin{align}
\label{2.11}
&
\displaystyle
\int_{t}^{x}(x-s)^{P-I}(s-t)^{Q-I}\ {}_{1}F_{2}\left(A;B,C;z(s-t)\right)ds
\nonumber
\\
&=
\displaystyle
\Gamma(P)\Gamma(Q)\Gamma^{-1}(P+Q)(x-t)^{P+Q-I} \ {}_{2}F_{3}\left(A,P;B,B+Q,C;z(x-t)\right),
\end{align}
\begin{align}
\label{2.12}
&
\displaystyle
\int_{x}^{y}(y-t)^{P-I}(t-x)^{Q-I}\ {}_{1}F_{2}\left(A;B,C;t-x\right) dt
\nonumber
\\
&=
\displaystyle
\Gamma(P)\Gamma(Q)\Gamma^{-1}(P+Q)(y-x)^{P+Q-I} \ {}_{2}F_{3}\left(A,Q;P+Q,B,C;y-x\right)
\end{align}
and
\begin{align}
\label{2.13}
&
\displaystyle
\int_{x}^{y}(y-t)^{P-I}(t-x)^{Q-I}\ {}_{1}F_{2}\left(A;B,C;x-t\right) dt
\nonumber
\\
&=
\displaystyle
\Gamma(P)\Gamma(Q)\Gamma^{-1}(P+Q)(y-x)^{P+Q-I} \ {}_{2}F_{3}\left(A,Q;P+Q,B,C;x-y\right).
\end{align}
\end{theorem}

\begin{proof}
Changing the variable $s$ to $u=\frac {s-t}{x-t}$ in the left hand side of (\ref{2.11})
and then using (\ref{1.7}) and  (\ref{1.8}), we obtain
\begin{eqnarray*}
\begin{split}
&
\displaystyle
\int_{t}^{x}(x-s)^{P-I}(s-t)^{Q+(k-1)I}ds=(x-t)^{P+Q+(k-1)I}\int_{0}^{1}(1-u)^{P-I}u^{Q+(k-1)I}du
\\
&=
\displaystyle
(x-t)^{P+Q+(k-1)I}B(P,Q+kI)=(x-t)^{P+Q+(k-1)I}\Gamma(P)\Gamma(Q+kI)\Gamma^{-1}(P+Q+kI)
\\
&=
\displaystyle
(x-t)^{P+Q+(k-1)I}\Gamma(P)\Gamma(Q)\Gamma^{-1}(P+Q)(Q)_{k}\left[(P+Q)_{k}\right]^{-1}.
\end{split}
\end{eqnarray*}
Similarly, changing the variable $t$ to $u=\frac {t-x}{y-x}$ in the left hand side of (\ref{2.12})
and then using (\ref{1.7}) and  (\ref{1.8}), we obtain
\begin{eqnarray*}
\begin{split}
&
\displaystyle
\int_{x}^{y}(y-t)^{P-I}(t-x)^{Q+(k-1)I}dt=(x-t)^{P+Q+(k-1)I}\int_{0}^{1}(1-u)^{P-I}u^{Q+(k-1)I}du
\\
&=
\displaystyle
(y-x)^{P+Q+(k-1)I}B(P,Q+kI)=(y-x)^{P+Q+(k-1)I}\Gamma(P)\Gamma(Q+kI)\Gamma^{-1}(P+Q+kI)
\\
&=
\displaystyle
(y-x)^{P+Q+(k-1)I}\Gamma(P)\Gamma(Q)\Gamma^{-1}(P+Q)(Q)_{k}\left[(P+Q)_{k}\right]^{-1}.
\end{split}
\end{eqnarray*}
Interpreting the right  hand side of (\ref{2.6}), in the view of  (\ref{1.2}) and the
concept of the Hadamard given by (\ref{1.3}), we arrive at the required result (\ref{2.5}).
\end{proof}

The proof of (\ref{2.13}) follows the same lines as that above.

\begin{theorem}
\label{thm 2.7}
Let  $A$, $B$, $C$, $A_{1}$ and $B_{1}$ be commuting positive stable matrices in $\Bbb{C}^{N\times N}$.
Then, for $Re(p)\geq Re(\lambda)$, we have the following integral representations
\begin{align}
\label{2.16}
&
\displaystyle
\int_{0}^{\infty}e^{-pt}t^{2A_{1}-1}\ {}_{1}F_{2}\left(A;B,C;\lambda^{2}t^{2}\right)dt
\nonumber
\\
&=
\displaystyle
\Gamma\left(2A_{1}\right)p^{-2A_{1}}\ {}_{3}F_{2}\left( A,\frac {1}{2}A_{1},\frac {1}{2}\left(A_{1}+I\right);B,C;\frac {4\lambda^{2}}{p^{2}}\right),
\
p\neq0
\end{align}
and
\begin{align}
\label{2.17}
&
\displaystyle
\int_{0}^{x}s^{B_{1}-1}(x-s)^{A_{1}-1}\ {}_{1}F_{2}\left(A;B,C;z(x-s)\right)ds
\nonumber
\\
&=
\displaystyle
\Gamma\left(A_{1}\right)\Gamma\left(B_{1}\right)\Gamma^{-1}\left(A_{1}+B_{1}\right)
(x-t)^{A_{1}+B_{1}-I} \ {}_{2}F_{3}\left(A,A_{1};B,B+B_{1},C;zx\right).
\end{align}
\end{theorem}

\begin{proof}
By taking $u=pt$ in the left  hand side of (\ref{2.16}) and then using (\ref{1.6}), we obtain
\begin{align*}
\displaystyle
\int_{0}^{\infty}e^{-pt}t^{2A_{1}+(2k-1)I}dt
&=
\displaystyle
\int_{0}^{\infty}e^{-u}\left(\frac {u}{p}\right)^{2A_{1}+(2k-1)I}\frac {1}{p}du
\\
&=
\displaystyle
\frac {1}{p^{2A_{1}+2kI}}\int_{0}^{\infty}e^{-u}u^{2A_{1}+(2k-1)I}du
\\
&=
\displaystyle
\frac {\Gamma\left(2A_{1}+2kI\right)}{p^{2A_{1}+2kI}},
\
p\neq0,
\end{align*}
\begin{eqnarray*}
\displaystyle
\frac {\Gamma\left(2A_{1}+2kI\right)}{p^{2A_{1}+2kI}}=\frac {2^{2A_{1}+2kI}
\Gamma\left(A_{1}+kI\right)\Gamma\left(A_{1}+(k+\frac {1}{2})I\right)}{\sqrt{\pi}p^{2A_{1}+2kI}},
\end{eqnarray*}
\begin{eqnarray*}
\displaystyle
\Gamma\left(2A_{1}\right)=\frac {2^{2A_{1}}\Gamma\left(A_{1}\right)\Gamma\left(A_{1}+\frac {1}{2}I\right)}{\sqrt{\pi}}
\end{eqnarray*}
and
\begin{align*}
&
\displaystyle
\sum_{k=0}^{\infty}\frac {1}{k!}(A)_{k}\left[(B)_{k}\right]^{-1}
\left[(C)_{k}\right]^{-1}\lambda^{2k}\int_{0}^{\infty}e^{-pt}t^{2A_{1}+(2k-1)I}dt
\\
&=
\displaystyle
\sum_{k=0}^{\infty}\frac {1}{k!}(A)_{k}\left[(B)_{k}\right]^{-1}\left[(C)_{k}\right]^{-1}
\lambda^{2k}\frac {2^{2A_{1}+2kI}
\Gamma\left(A_{1}+kI\right)
\Gamma\left(A_{1}+\left(k+\frac {1}{2}\right)I\right)}{\sqrt{\pi}p^{2A_{1}+2kI}}
\\
&=
\displaystyle
\Gamma\left(A_{1}\right)\Gamma\left(A_{1}+\frac {1}{2}\right)
\sum_{k=0}^{\infty}
\frac {1}{k!}(A)_{k}\left(A_{1}\right)_{k}
\left(A_{1}+\frac {1}{2}\right)_{k}
\left[(B)_{k}\right]^{-1}
\left[(C)_{k}\right]^{-1}\frac {2^{2A_{1}+2k}\lambda^{2k}}{\sqrt{\pi}p^{2A_{1}+2kI}}
\\
&=
\displaystyle
\frac {2^{2A_{1}}}{\sqrt{\pi}p^{2A_{1}}}
\Gamma\left(A_{1}\right)
\Gamma\left(A_{1}+\frac {1}{2}\right)
\sum_{k=0}^{\infty}
\frac {1}{k!}(A)_{k}\left(A_{1}\right)_{k}
\left(A_{1}+\frac {1}{2}I\right)_{k}
\left[(B)_{k}\right]^{-1}
\left[(C)_{k}\right]^{-1}\frac {2^{2k}\lambda^{2k}}{p^{2k}}
\\
&=
\displaystyle
\frac {\Gamma(2A_{1})}{p^{2A_{1}}}
\sum_{k=0}^{\infty}
\frac {1}{k!}(A)_{k}\left(A_{1}\right)_{k}
\left(A_{1}+\frac {1}{2}\right)_{k}
\left[(B)_{k}\right]^{-1}
\left[(C)_{k}\right]^{-1}\left(\frac {2\lambda}{p}\right)^{2k}.
\end{align*}
By using the concept of the  generalized hypergeometric matrix function on the left  hand side of the above relation,
we obtain the right hand side of (\ref{2.16}).
Similarly, we can prove (\ref{2.17}).
\end{proof}

\begin{theorem}
\label{thm 3.2}
Let $A$, $B$, $C$, $P$ and $Q$ be commutative matrices in $\Bbb{C}^{N\times N}$ satisfying  (\ref{2.1}).
Then, the following holds
\begin{equation*}
\displaystyle
\int_{0}^{x}t^{P-I}(x-t)^{Q-I}\ {}_{1}F_{2}(A;B,C;t)dt=x^{P+Q-I}\mathbf{B}(P,Q)\ {}_{2}F_{3}(A,P;B,C,P+Q;x).
\end{equation*}
\end{theorem}

\begin{proof}
Similar to the proof of Theorem \ref{thm 2.7}.
\end{proof}

\begin{theorem}
\label{thm 2.9}
Let $A$, $B$ and $C$ be matrices in $\Bbb{C}^{N\times N}$ satisfying (\ref{2.1}) and all matrices are commutative.
Then,
\begin{eqnarray}
\label{2.19}
\begin{split}
\displaystyle
(\theta\ I+A)\ {}_{1}F_{2}(A;B,C;z)=A \ {}_{1}F_{2}(A+I;B,C;z),
\\
\displaystyle
(\theta\ I+B-I)\ {}_{1}F_{2}(A;B,C;z)=(B-I)\ {}_{1}F_{2}(A;B-I,C;z),
\\
\displaystyle
(\theta\ I+C-I)\ {}_{1}F_{2}(A;B,C;z)=(C-I)\ {}_{1}F_{2}(A;B,C-I;z),
\end{split}
\end{eqnarray}
where $\theta$ is the Euler differential operator $z\frac {d}{z}$ and
\begin{eqnarray*}
\begin{split}
\displaystyle
(A-B+I)\ {}_{1}F_{2}(A;B,C;z)=A\ {}_{1}F_{2}(A+I;B,C;z)-(B-I)\ {}_{1}F_{2}(A;B-I,C;z),
\\
\displaystyle
(A-C+I)\ {}_{1}F_{2}(A;B,C;z)=A\ {}_{1}F_{2}(A+I;B,C;z)-(C-I)\ {}_{1}F_{2}(A;B,C-I;z),
\\
\displaystyle
(B-C)\ {}_{1}F_{2}(A;B,C;z)=(B-I)\ {}_{1}F_{2}(A;B-I,C;z)-(C-I)\ {}_{1}F_{2}(A;B,C-I;z).
\end{split}
\end{eqnarray*}
\end{theorem}

\begin{theorem}
\label{thm 2.9}
Let $A$, $B$ and $C$ be commutative matrices in $\Bbb{C}^{N\times N}$ satisfying  (\ref{2.1}).
Then ${}_{1}F_{2}$ is a solution of
\begin{eqnarray}
\displaystyle
\left[  \theta\ (\theta\ I+B-I)(\theta\ I+C-I)-z(\theta\ I+A)\right]\ {}_{1}F_{2}=\mathbf{0}.
\label{2.20}
\end{eqnarray}
\end{theorem}

\begin{proof}
From (\ref{2.2}) and (\ref{2.19}), if we take in consideration the differential operator $\theta=z\frac {d}{dz}$, $\theta{z^{k}}=k{z^{k}}$.
Thus,
\begin{eqnarray*}
\displaystyle
\theta\ (\theta\ I+B-I)(\theta\ I+C-I)\ {}_{1}F_{2}=
\sum_{k=1}^{\infty}\frac {z^{k}}{(k-1)!}(A)_{k}
\left[(B)_{k-1}\right]^{-1}
\left[(C)_{k-1}\right]^{-1}.
\end{eqnarray*}
Now replace $k$ by $k+1$, we have
\begin{align*}
\displaystyle
\theta\ (\theta\ I+B-I)(\theta\ I+C-I)\ {}_{1}F_{2}
&=
\displaystyle
\sum_{k=0}^{\infty}\frac {z^{k+1}}{k!}(A)_{k+1}
\left[(B)_{k}\right]^{-1}
\left[(C)_{k}\right]^{-1}
\\
&=
\displaystyle
\sum_{k=0}^{\infty}\frac {z^{k+1}}{k!}(A+kI)(A)_{k}
\left[(B)_{k}\right]^{-1}
\left[(C)_{k}\right]^{-1}
\\
&=
\displaystyle
z(\theta\ I+A)\ {}_{1}F_{2}.
\end{align*}
The proof is complete.
\end{proof}

\begin{theorem}
\label{thm 2.9}
Let $A$, $B$ and $C$ be commutative matrices in $\Bbb{C}^{N\times N}$ satisfying (\ref{2.1}).
Then ${}_{1}F_{2}$ is a solution of
\begin{eqnarray*}
\displaystyle
(\theta\ I+B)(\theta\ I+C)(\theta\ I+I)\left( \frac {{}_{1}F_{2}}{z}\right) - (\theta\ I+A)\ {}_{1}F_{2}=\mathbf{0}.
\end{eqnarray*}
\end{theorem}

\begin{theorem}
\label{thm 2.10}
Let $A$, $B$ and $C$ be commutative matrices in $\Bbb{C}^{N\times N}$ satisfying  (\ref{2.1}).
Then ${}_{1}F_{2}$ is a solution of
\begin{eqnarray*}
\displaystyle
z^{2}\frac {d^{3}}{dz^{3}}\ {}_{1}F_{2}+(B+C+I)z\frac {d^{2}}{dz^{2}}\ {}_{1}F_{2}+(BC-z I)\frac {d}{dz}\ {}_{1}F_{2}-A\ {}_{1}F_{2}=\mathbf{0}.
\end{eqnarray*}
\end{theorem}

\begin{theorem}
\label{thm 2.11}
Let $A$, $B$ and $C$ be commutative matrices in $\Bbb{C}^{N\times N}$ satisfying  (\ref{2.1}).
Then we have
\begin{align}
&
\displaystyle
z^{2}\frac {d^{3}}{dz^{3}}\ \Psi(z)+(3\alpha I-B-C+I)z\frac {d^{2}}{dz^{2}}\ \Psi(z)+\bigg[3\alpha(\alpha-1)I+2\alpha(B+C+I)
\nonumber
\\
&
\displaystyle
+(BC-z I)\bigg]\frac {d}{dz}\ \Psi(z) -\bigg[\alpha I+A-\frac {\alpha}{z}\left((\alpha-1)I-B\right)
\left((\alpha-1)I+C\right)\bigg] \ \Psi(z)=\mathbf{0}.
\label{2.23}
\end{align}
\end{theorem}

\begin{proof}
Let $Y(z)=\sum_{k=0}^{\infty}c_{k}z^{k+\alpha}$.
Then,
\begin{eqnarray*}
\displaystyle
Y'(z)=\sum_{k=0}^{\infty}(k+\alpha)C_{k}z^{k+\alpha-1},
\
\displaystyle
Y''(z)=\sum_{k=0}^{\infty}(k+\alpha)(k+\alpha-1)C_{k}z^{k+\alpha-2}.
\end{eqnarray*}
Substituting in (\ref{2.17}), rearranging terms and summations, we have three cases to consider:
$k=0$, $C_{0}\neq \mathbf{0}$, $\alpha=\mathbf{0}$, $\alpha=I-B$ or $\alpha=I-C$,  ${}_{1}F_{2}=z^{\alpha} \Psi(z)$.

Case (I), $\alpha=\mathbf{0}$
\begin{eqnarray*}
\displaystyle
Y_{1}(z)=\ {}_{1}F_{2}(A;B,C;z).
\end{eqnarray*}

Case (II), $\alpha=I-B$
\begin{eqnarray*}
\displaystyle
Y_{2}(z)=z^{I-B}\ {}_{1}F_{2}(I-B+A;2I-B,I-C+B;z).
\end{eqnarray*}

Case (III), $\alpha=I-C$
\begin{eqnarray*}
\displaystyle
Y_{3}(z)=z^{I-C}\  {}_{1}F_{2}(I-C+A;I+B-C,2I-C;z).
\end{eqnarray*}

Hence, the general solution is
\begin{align*}
\displaystyle
Y(z)
&=
\displaystyle
a_{1}Y_{1}(z)+a_{2}Y_{2}(z)+a_{3}Y_{3}(z)=a_{1}\ {}_{1}F_{2}(A;B,C;z)
\\
&
\displaystyle
+a_{2}z^{I-B}\ {}_{1}F_{2}(I-B+A;2I-B,I-C+B;z)+a_{3}z^{I-C}\ {}_{1}F_{2}(I-C+A;I+B-C,2I-C;z),
\end{align*}
where $a_{1}$, $a_{2}$ and $a_{3}$ are arbitrary constants.
\end{proof}

\begin{theorem}
\label{thm 2.12}
Let $A$, $B$ and $C $ be the positive stable matrices in $ \Bbb{C}^{N \times N}$ satisfying  (\ref{2.1}).
Then, for $ |z| < 1 $,  the following hypergeometric matrix functions ${}_{1}F_{2}$ hold true
\begin{eqnarray}
\label{2.28}
\displaystyle
{}_{1}F_{2}(A;B,C;z)=\sum_{k=0}^{\infty}\frac {(-1)^{k}(B-A)_{k}}{k!}
\left[(B)_{k}\right]^{-1}\left[(C)_{k}\right]^{-1}z^{k} \ {}_{0}F_{1}(-;C+k I;z),
\end{eqnarray}
where $\Re(B)>\Re(A)>0$,
\begin{eqnarray}
\label{2.29}
\displaystyle
{}_{1}F_{2}(A;B,C;z)=\sum_{k=0}^{\infty}\frac {(-1)^{k}(C-A)_{k}}{k!}
\left[(B)_{k}\right]^{-1}
\left[(C)_{k}\right]^{-1}z^{k}
\ {}_{0}F_{1}(-;B+k I;z),
\end{eqnarray}
where $\Re(C)>\Re(A)>0$,
\begin{eqnarray}
\label{2.30}
\displaystyle
{}_{0}F_{1}(-;C;z)=\sum_{k=0}^{\infty}
\frac {(B-A)_{k}}{k!}[(B)_{k}]^{-1}\left[(C)_{k}\right]^{-1}z^{k}
\ {}_{1}F_{2}(A;B+k I,C+k I;z),
\end{eqnarray}
where $\Re(B)>\Re(A)>0$, and
\begin{eqnarray}
\label{2.31}
\displaystyle
{}_{0}F_{1}(-;B;z)=\sum_{k=0}^{\infty}
\frac {(C-A)_{k}}{k!}\left[(B)_{k}\right]^{-1}\left[(C)_{k}\right]^{-1}z^{k}\ {}_{1}F_{2}(A;B+k I,C+k I;z),
\end{eqnarray}
where $\Re(C)>\Re(A)>0$.
\end{theorem}

\begin{theorem}
\label{thm 2.14}
The following  infinite summations for ${}_{1}F_{2}$ hold true
\begin{eqnarray}
\label{2.311}
\displaystyle
\sum_{k=0}^{\infty}\frac {(A)_{k}}{k!}\ {}_{1}F_{2}(A+kI;B,C;z)t^{k}=(1-t)^{-A}\ {}_{1}F_{2}(A;B,C;\frac {z}{1-t}),
\
|t|<1,
\end{eqnarray}
\begin{eqnarray}
\label{2.3111}
\displaystyle
\sum_{k=0}^{\infty}\frac {(A)_{k}\left[(B)_{k}\right]^{-1}
\left[(C)_{k}\right]^{-1}}{k!}\ {}_{1}F_{2}(A+kI;B+kI,C+kI;z)t^{k}=\ {}_{1}F_{2}(A;B,C;z+t)
\end{eqnarray}
and
\begin{eqnarray}
\label{2.32}
\displaystyle
\sum_{k=0}^{\infty}\frac {(E)_{k}\left[(B)_{k}\right]^{-1}\left[(C)_{k}\right]^{-1}}{k!}
\ {}_{1}F_{2}(A;B+kI,C+kI;z)z^{k}=(1-t)^{-A}\ {}_{1}F_{2}(A+E;B,C;z),
\
|t|<1.
\end{eqnarray}
\end{theorem}

\section{Fractional integrals and derivatives of ${}_{1}F_{2}$}

In this section, we study certain integral transforms and fractional integral formulas for  ${}_{1}F_{2}$.
The  beta transform of $f(x)$ for all $t\geq0$ is defined as (see \cite{ba1, ba3, ba4, ba5, ba6})
\begin{eqnarray*}
\displaystyle
B\left\{ f(x);A,B\right\}=\int_{0}^{1}t^{A-I}(1-t)^{B-I}f(t)dt.
\end{eqnarray*}

The Laplace transform of  $f(z)$ is defined  by
\begin{equation*}
\displaystyle
\mathbb{L}\left\{f(z):s\right\}=\int_{0}^{\infty}e^{-sz}f(z)dz,
\
\displaystyle
\Re(s)>0,
\end{equation*}
where $f\in L(0,R)$ for every $R>0$ and $f(t)=O\left(e^{at}\right)$, $t\rightarrow \infty$.

The Erdelyi Kober type fractional integral operators are defined as
\begin{eqnarray*}
\displaystyle
\mathbb{E}^{\alpha,\eta}_{0,x}f(x)=\frac {x^{-\eta-\alpha}}{\Gamma(\alpha)}\int_{0}^{x}(x-t)^{\alpha-1}t^{\eta}f(t)dt
\end{eqnarray*}
and
\begin{eqnarray*}
\displaystyle
\mathbb{K}^{-\alpha,\eta}_{x,\infty}f(x)=\frac {x^{\eta}}{\Gamma(\alpha)}\int_{x}^{\infty}(t-x)^{\alpha-1}t^{-\eta-\alpha}f(t)dt,
\end{eqnarray*}
where  $f(t)$ is so constrained that the defining integrals  converge.

The fractional integral and derivative of Riemann-Liouville of order $\mu $ are defined as follows (see $\cite{lren}$)
\begin{equation}
\label{3.5}
\displaystyle
\mathbb{I}^{\mu}\left\{f(x)\right\}=\frac {1}{\Gamma(\mu)}\int_{0}^{x}(x-t)^{\mu-1}f(t)dt,
\
\displaystyle
Re(\mu)>0
\end{equation}
and
\begin{equation*}
\displaystyle
\mathbb{D}^{\mu}= D^{n} \left(\mathbb{I}^{{n-\mu}}f(x)\right),
\
D=\frac {d}{dx}
\end{equation*}
for $x > 0$.

The classical Riemann–Liouville fractional derivative of order $\mu$ is defined by
\begin{equation}
\label{3.7}
\displaystyle
D_{z}^{\mu}f(z)=\frac {1}{\Gamma(-\mu)}\int_{0}^{z}\frac {f(t)}{(z-t)^{\mu+1}}dt,
\
\displaystyle
\Re(z)<0,
\end{equation}
where the integration path is a line from $0$ to $z$ in the complex $t$-plane.

Let $f(x)\in L(b,c)$, $\alpha\in\Bbb{C}$ and $Re(\alpha)>0$,
then the left sided and right sided operators
of Riemann-Liouville fractional integral of order $\alpha$ are defined as
\begin{equation}
\displaystyle
{}_{b}\mathbb{I}^{\alpha}_{x}\{f(x)\}=\frac {1}{\Gamma(\alpha)}\int_{b}^{x}(x-t)^{\alpha-1}f(t)dt,
\
\displaystyle
x>b
\label{3.8}
\end{equation}
and
\begin{equation}
\displaystyle
{}_{x}\mathbb{I}^{\alpha}_{c}\{f(x)\}=\frac {1}{\Gamma(\alpha)}\int_{x}^{c}(t-x)^{\alpha-1}f(t)dt,
\
\displaystyle
x<c,
\label{3.9}
\end{equation}
respectively.

Let $f(x)\in L(b,c)$, $\alpha\in\Bbb{C}$, $Re(\alpha)\geq0$ and $n=\left[Re(\alpha)\right]+1$.
The left sided and right sided operators of Riemann-Liouville fractional derivative of order $\alpha$ are defined by
\begin{equation}
\label{3.10}
\displaystyle
{}_{b}D_{x}^{\alpha}\left\{f(x)\right\}=\frac {1}{\Gamma(n-\alpha)}\left(\frac {\partial}{\partial x}\right)^{n}
\int_{b}^{x}\frac {f(t)}{(x-t)^{\alpha-n+1}}dt,
\
\displaystyle
x>b
\end{equation}
and
\begin{equation}
\label{3.11}
\displaystyle
{}_{x}D_{c}^{\alpha}\left\{f(x)\right\}=\frac {(-1)^{n}}{\Gamma(n-\alpha)}\left(\frac {\partial}{\partial\ x}\right)^{n}
\int_{x}^{c}\frac {f(t)}{(t-x)^{\alpha-n+1}}dt,
\
\displaystyle
x<c,
\end{equation}
respectively.

Let $f(x)\in L(b,c)$, $\alpha\in\Bbb{C}$, $Re(\alpha)\geq0$ and $n=\left[Re(\alpha)\right]+1$,
then the Weyl fractional derivative of $f(x)$ of order $\alpha$, denoted by
${}_{x}D_{\infty}^{\alpha}$, is defined by
\begin{equation}
\label{3.12}
\displaystyle
{}_{x}D_{\infty}^{\alpha}\{f(x)\}=
\frac {(-1)^{n}}{\Gamma(n-\alpha)}\left(
\frac {\partial}{\partial x}\right)^{n}\int_{x}^{\infty}\frac {f(t)}{(t-x)^{\alpha-n+1}}dt.
\end{equation}

The Weyl integral of $f(x)$ of order $\alpha$, denoted by ${}_{x}W_{\infty}^{\alpha}$, is defined by
\begin{equation}
\label{3.13}
\displaystyle
{}_{x}W_{\infty}^{\alpha}\left\{f(x)\right\}=\frac {1}{\Gamma(\alpha)}\int_{x}^{\infty}(t-x)^{\alpha-1}f(t)dt,
\
\displaystyle
-\infty<x<\infty.
\end{equation}

\begin{theorem}
\label{thm 3.1}
The following beta and Laplace transforms for ${}_{1}F_{2}$ in (\ref{2.2}) hold true
\begin{eqnarray}
\label{3.14}
\displaystyle
B\left\{ {}_{1}F_{2}(A;B,C;z);P,Q\right\}=B(P,Q)\ {}_{2}F_{3}(A,P;B,C,P+Q;z)
\end{eqnarray}
and
\begin{eqnarray}
\label{3.15}
\displaystyle
L\left\{z^{P-I}\ {}_{1}F_{2}(A;B,C;z);P,Q\right\}=\Gamma(P)s^{-P}\ {}_{2}F_{2}\left(A,P;B,C;\frac {z}{s}\right).
\end{eqnarray}
\end{theorem}

\begin{proof}
Using the results
\begin{eqnarray*}
\displaystyle
\int_{0}^{1}t^{P+(k-1)I}(1-t)^{Q-I}dt=\mathbf{B}(P+kI,Q)=\Gamma(P)\Gamma(Q)\Gamma^{-1}(P+Q)(P)_{k}\left[(P+Q)_{k}\right]^{-1}
\end{eqnarray*}
and
\begin{eqnarray*}
\displaystyle
\int_{0}^{\infty}t^{P+(k-1)I}e^{-st}dt=\Gamma(P+kI)s^{-(P+kI)}=\Gamma(P)s^{-(P+kI)}(P)_{k},
\end{eqnarray*}
we obtain (\ref{3.14})-(\ref{3.15}).
\end{proof}

\begin{theorem}
\label{thm 3.2}
The following  Erdelyi Kober type fractional integral operators for ${}_{1}F_{2}$ in (\ref{2.2}) hold true
\begin{align}
\label{3.16}
&
\displaystyle
\mathbb{E}^{\alpha,\eta}_{0,x}\ {}_{1}F_{2}(A;B,C;x)
\nonumber
\\
&=
\displaystyle
\frac {x^{-\eta-\alpha}}{\Gamma(\alpha)}
\int_{0}^{x}(x-t)^{\alpha-1}t^{\eta}\ {}_{1}F_{2}(A;B,C;t)dt
\nonumber
\\
&=
\displaystyle
\frac {\Gamma(\eta)}{\Gamma(\alpha+\eta)}\ {}_{2}F_{3}\left(A,\eta I;B,C,(\alpha+\eta)I;x\right)
\end{align}
and
\begin{align}
\label{3.17}
&
\displaystyle
\mathbb{K}^{-\alpha,\eta}_{x,\infty}\ {}_{1}F_{2}(A;B,C;x)
\nonumber
\\
&=
\displaystyle
\frac {x^{\eta}}{\Gamma(\alpha)}\int_{x}^{\infty}(x-t)^{\alpha-1}t^{-\eta-\alpha}\ {}_{1}F_{2}(A;B,C;t)dt
\nonumber
\\
&=
\displaystyle
\frac {\Gamma(\eta)}{\Gamma(\alpha+\eta)}\ {}_{2}F_{3}\left(A,\eta I;B,C,(\alpha+\eta)I;x\right).
\end{align}
\end{theorem}

\begin{proof}
Put $t=xu$ and $dt=xdu$, we have
\begin{eqnarray}
\displaystyle
\mathbb{E}^{\alpha,\eta}_{0,x}x^{k}=\frac {x^{-\eta-\alpha}}{\Gamma(\alpha)}
\int_{0}^{x}(x-t)^{\alpha-1}t^{\eta+k}dt =
\frac {\Gamma(\eta+k)}{\Gamma(\alpha+\eta+k)}x^{k} =
\frac {\Gamma(\eta)(\eta)_{k}}{\Gamma(\alpha+\eta)(\alpha+\eta)_{k}}x^{k}.
\label{3.18}
\end{eqnarray}
By using  (\ref{3.18}) and (\ref{2.2}), we obtain (\ref{3.16}).

Put $t=\frac {x}{u}$ and $dt=-\frac {x}{u^{2}}du$, we obtain
\begin{eqnarray*}
\displaystyle
\mathbb{K}^{-\alpha,\eta}_{x,\infty}x^{k} =
\frac {x^{\eta}}{\Gamma(\alpha)}\int_{x}^{\infty}(t-x)^{\alpha-1}t^{k-\eta-\alpha}dt =
\frac {\Gamma(\eta+k)}{\Gamma(\alpha+\eta+k)}x^{k} =
\frac {\Gamma(\eta)(\eta)_{k}}{\Gamma(\alpha+\eta)(\alpha+\eta)_{k}}x^{k}.
\end{eqnarray*}
By using (\ref{2.2}), we obtain (\ref{3.17}).
\end{proof}

\begin{theorem}
\label{thm 3.3}
The following  fractional integral operators for ${}_{1}F_{2}$ in (\ref{2.2}) hold true
\begin{equation}
\label{3.20}
\displaystyle
\mathbb{I}^{\mu}\left\{ \ {}_{1}F_{2}(A;B,C;x)\right\}=\frac {x^{\mu}}{\Gamma(\mu+1)}\ {}_{2}F_{3}\left(A,I;B,C,(\mu+1)I;x\right),
\end{equation}
\begin{equation}
\label{3.21}
\displaystyle
{}_{b}\mathbb{I}^{\alpha}_{x}\left\{f(x)\right\} =\frac {(x-b)^{\alpha}}{\Gamma(\alpha+1)}
\ {}_{2}F_{3}\left(A,I;B,C,(\alpha+1)I;x-b\right),
\end{equation}
\begin{equation}
\label{3.22}
\displaystyle
{}_{x}\mathbb{I}^{\alpha}_{a}\left\{ \ {}_{1}F_{2}(A;B,C;x)\right\} =
\frac {(a-x)^{\alpha}}{\Gamma(\alpha+1)}\  {}_{2}F_{3}\left(A,I;B,C,(\alpha+1)I;a-x\right)
\end{equation}
and
\begin{equation}
\displaystyle
{}_{x}W_{\infty}^{\alpha}\left\{\ {}_{1}F_{2}(A;B,C;x)\right\}=x^{\alpha}\ {}_{1}F_{3}(A;B,C,\alpha I;x),
\label{3.23}
\end{equation}
where $\alpha\in\Bbb{C}$ and $Re(\alpha)>0$.
\end{theorem}

\begin{proof}
Put $t=xu$  and $dt=xdu$, we have
\begin{align}
&
\displaystyle
\frac {1}{\Gamma(\mu)}\int_{0}^{x}(x-t)^{\mu-1}t^{k}dt
\nonumber
\\
&=
\displaystyle
\frac {1}{\Gamma(\mu)}x^{\mu+k-1+1}\int_{0}^{1}(1-u)^{\mu-1}u^{k}du=\frac {1}{\Gamma(\mu)}x^{\mu+k-1+1}B(\mu,k+1)
\nonumber
\\
&=
\displaystyle
\frac {1}{\Gamma(\mu)}x^{\mu+k}\frac {\Gamma(\mu)\Gamma(k+1)}{\Gamma(\mu+k+1)}=x^{\mu+k}\frac {(1)_{k}}{\Gamma(\mu+1)(\mu+1)){k}}.
\label{3.24}
\end{align}
Using (\ref{3.24}), (\ref{3.5}) and (\ref{2.2}), we obtain (\ref{3.20}).

Put $u=\frac {b-t}{b-x}$ and $dt=-(b-x)du$, we have
\begin{align}
&
\displaystyle
\frac {1}{\Gamma(\alpha)}\int_{b}^{x}(x-t)^{\alpha-1}(t-b)^{k}dt
\nonumber
\\
&=
\displaystyle
\frac {1}{\Gamma(\alpha)}\int_{0}^{1}(x-b+(b-x)u)^{\alpha-1}\left[(b-x)u\right]^{k}dt
\nonumber
\\
&=
\displaystyle
\frac {1}{\Gamma(\alpha)}(x-b)^{\alpha+k-1+1}\int_{0}^{1}(1-u)^{\alpha-1}u^{k}dt
\nonumber
\\
&=
\displaystyle
\frac {1}{\Gamma(\alpha)}(x-b)^{\alpha+k}B(\alpha,k+1)
\nonumber
\\
&=
\displaystyle
\frac {1}{\Gamma(\alpha)}(x-b)^{\alpha+k}\frac {\Gamma(\alpha)\Gamma(k+1)}{\Gamma(\alpha+k+1)}
\nonumber
\\
&=
\displaystyle
(x-b)^{\alpha+k}\frac {(1)_{k}}{\Gamma(\alpha+1)(\alpha+1)_{k}}.
\label{3.25}
\end{align}
Using (\ref{3.25}), (\ref{3.8}) and (\ref{2.2}), we obtain (\ref{3.21}).

Put $u=\frac {c-t}{c-x}$ and $dt=-(c-x)du$, we have
\begin{align}
&
\displaystyle
\frac {1}{\Gamma(\alpha)}\int_{x}^{c}(t-x)^{\alpha-1}(c-t)^{k}dt
\nonumber
\\
&=
\displaystyle
\frac {1}{\Gamma(\alpha)}(c-x)^{\alpha+k}\int_{0}^{1}(1-u)^{\alpha-1}u^{k}du
\nonumber
\\
&=
\displaystyle
\frac {1}{\Gamma(\alpha)}(c-x)^{\alpha+k}\frac {\Gamma(\alpha)\Gamma(k+1)}{\Gamma(\alpha+k+1)}
\nonumber
\\
&=
\displaystyle
(c-x)^{\alpha+k}\frac {(1)_{k}}{\Gamma(\alpha+1)(\alpha+1)_{k}}.
\label{3.26}
\end{align}
Using (\ref{3.26}), (\ref{3.9}) and (\ref{2.2}), we obtain (\ref{3.22}).

Put $t=\frac {x}{u}$  and $dt=-\frac {x}{u^{2}}du$, we have
\begin{equation}
\displaystyle
\mathbb{W}^{-\alpha}_{x,\infty}x^{k}
=\frac {1}{\Gamma(\alpha)}\int_{x}^{\infty}(t-x)^{\alpha-1}t^{k}dt
=\frac {\Gamma(k)}{\Gamma(\alpha+k)}x^{\alpha+k}=\frac {1}{(\alpha)_{k}}x^{\alpha+k}.
\label{3.27}
\end{equation}
Using (\ref{3.27}), (\ref{3.13}) and (\ref{2.2}), we obtain (\ref{3.23}).
\end{proof}

\begin{theorem}
\label{thm 3.4}
The following  fractional derivatives operators for ${}_{1}F_{2}$ in (\ref{2.2}) hold true
\begin{equation}
\displaystyle
{}_{b}D_{x}^{\alpha}\left\{ {}_{1}F_{2}(A;B,C;x-b)\right\}=\frac {(x-b)^{-\alpha}}{\Gamma(1-\alpha)}\ {}_{2}F_{3}\left(A,I;B,C,(1-\alpha)I;x-b\right),
\label{3.28}
\end{equation}
\begin{equation}
\displaystyle
{}_{x}D_{c}^{\alpha}\left\{{}_{1}F_{2}(A;B,C;c-x)\right\}=\frac {(c-x)^{-\alpha}}{\Gamma(1-\alpha)}
\ {}_{2}F_{3}\left(A,I;B,C,(1-\alpha)I;c-x\right),
\label{3.29}
\end{equation}
\begin{equation}
\displaystyle
{}_{x}D_{\infty}^{\alpha}\left\{ {}_{1}F_{2}(A;B,C;x)\right\} =
\frac {(-1)^{\alpha}x^{-\alpha}}{\Gamma(1-\alpha)}\ {}_{2}F_{3}\left(A,I;B,C,(1-\alpha)I;x\right)
\label{3.30}
\end{equation}
and
\begin{equation}
\displaystyle
D_{z}^{\mu}\ {}_{1}F_{2}(A;B,C;x)=\frac {x^{-\mu}}{\Gamma(1-\mu)}\ {}_{2}F_{3}\left(A,I;B,C,(1-\mu)I;x\right).
\label{3.31}
\end{equation}
\end{theorem}

\begin{proof}
Put $u=\frac {b-t}{b-x}$ and $dt=-(b-x)du$,  we have
\begin{align}
&
\displaystyle
\frac {1}{\Gamma(n-\alpha)}\int_{b}^{x}\frac {(t-b)^{k}}{(x-t)^{\alpha-n+1}}dt
\nonumber
\\
&=
\displaystyle
\frac {1}{\Gamma(n-\alpha)}\int_{0}^{1}(x-b+(b-x)u)^{n-\alpha-1}\left[(x-b)u\right]^{k}dt
\nonumber
\\
&=
\displaystyle
\frac {1}{\Gamma(n-\alpha)}(x-b)^{n-\alpha+k-1+1}\int_{0}^{1}(1-u)^{n-\alpha-1}u^{k}ut
\nonumber
\\
&=
\displaystyle
\frac {1}{\Gamma(\alpha)}(x-b)^{k+n-\alpha}B(n-\alpha,k+1)
\nonumber
\\
&=
\displaystyle
\frac {1}{\Gamma(n-\alpha)}(x-b)^{k+n-\alpha}\frac {\Gamma(n-\alpha)\Gamma(k+1)}{\Gamma(n-\alpha+k+1)}.
\label{3.32}
\end{align}
Using (\ref{3.32}), (\ref{3.10}) and (\ref{2.2}), we obtain (\ref{3.28}).

Put $u=\frac {c-t}{c-x}$  and $dt=-(c-x)du$,  we have
\begin{align}
&
\displaystyle
\frac {(-1)^{n}}{\Gamma(n-\alpha)}\int_{x}^{c}\frac {(c-t)^{k}}{(t-x)^{\alpha-n+1}}dt
\nonumber
\\
&=
\displaystyle
\frac {1}{\Gamma(n-\alpha)}(c-x)^{n-\alpha+k}\int_{0}^{1}(1-u)^{n-\alpha-1}u^{k}du
\nonumber
\\
&=
\displaystyle
\frac {(-1)^{n}}{\Gamma(n-\alpha)}(c-x)^{n-\alpha+k}B(n-\alpha,k+1)
\nonumber
\\
&=
\displaystyle
\frac {(-1)^{n}}{\Gamma(n-\alpha)}(c-x)^{n-\alpha+k}\frac {\Gamma(n-\alpha)\Gamma(k+1)}{\Gamma(n-\alpha+k+1)}
\nonumber
\\
&=
\displaystyle
(-1)^{n}(c-x)^{n-\alpha+k}\frac {\Gamma(k+1)}{\Gamma(n-\alpha+k+1)}.
\label{3.33}
\end{align}
Using (\ref{3.33}), (\ref{3.11}) and (\ref{2.2}), we obtain (\ref{3.29}).

Put $t=\frac {x}{u}$  and $dt=-\frac {x}{u^{2}}du$,  we have
\begin{align}
&
\displaystyle
\mathbb{W}^{-\alpha}_{x,\infty}x^{k}=\frac {(-1)^{n}}{\Gamma(n-\alpha)}D^{n}\int_{x}^{\infty}\frac {t^{k}}{(t-x)^{\alpha-n+1}}dt
\nonumber
\\
&=
\displaystyle
D^{n}\frac {(-1)^{n}}{\Gamma(n-\alpha)}\int_{1}^{0}
\frac {\left(\frac {x}{u}\right)^{k}}{\left(\frac {x}{u}-x\right)^{\alpha-n+1}}\left(-\frac {x}{u^{2}}du\right)
\nonumber
\\
&=
\displaystyle
D^{n}x^{n-\alpha+k}\frac {(-1)^{n}}{\Gamma(n-\alpha)}
\int_{0}^{1}u^{\alpha-n-k-1}(1-u)^{n-\alpha-1}du)
\nonumber
\\
&=
\displaystyle
D^{n}x^{n-\alpha+k}\frac {(-1)^{n}}{\Gamma(n-\alpha)}B(\alpha-n-k,n-\alpha)
\nonumber
\\
&=
\displaystyle
D^{n}x^{n-\alpha+k}\frac {(-1)^{n}}{\Gamma(n-\alpha)}
\frac {\Gamma(\alpha-n-k)\Gamma(n-\alpha)}{\Gamma(-k)}
\nonumber
\\
&=
\displaystyle
\frac {\Gamma(n-\alpha+k+1)}{\Gamma(k-\alpha+1)}x^{k-\alpha}\frac {(-1)^{n}(-1)^{n+k-\alpha}\Gamma(k+1)}{(-1)^{k}\Gamma(n+k-\alpha+1)}
\nonumber
\\
&=
\displaystyle
\frac {(-1)^{\alpha}\Gamma(k+1)}{\Gamma(k-\alpha+1)}x^{k-\alpha}
\nonumber
\\
&=
\displaystyle
\frac {(-1)^{\alpha}(1)_{k}}{\Gamma(1-\alpha)(1-\alpha)_{k}}x^{k-\alpha}.
\label{3.34}
\end{align}
Using (\ref{3.34}), (\ref{3.12}) and (\ref{2.2}), we obtain (\ref{3.30}).

Put $t=xu$  and $dt=xdu$, we have
\begin{align}
&
\displaystyle
\frac {1}{\Gamma(-\mu)}\int_{0}^{x}\frac {t^{k}}{(x-t)^{\mu+1}}dt
\nonumber
\\
&=
\displaystyle
\frac {1}{\Gamma(-\mu)}x^{-\mu+k-1+1}\int_{0}^{1}(1-u)^{-\mu-1}u^{k}du
\nonumber
\\
&=
\displaystyle
\frac {1}{\Gamma(-\mu)}x^{k-\mu}B(-\mu,k+1)
\nonumber
\\
&=
\displaystyle
\frac {1}{\Gamma(-\mu)}x^{k-\mu}\frac {\Gamma(-\mu)\Gamma(k+1)}{\Gamma(-\mu+k+1)}
\nonumber
\\
&=
\displaystyle
x^{k-\mu}\frac {(1)_{k}}{\Gamma(1-\mu)(1-\mu)){k}}.
\label{3.35}
\end{align}
Using (\ref{3.35}), (\ref{3.7}) and (\ref{2.2}), we obtain (\ref{3.31}).
\end{proof}

\section{Applications}

\subsection{Bateman matrix polynomials $\mathbf{B}_{n}^{A,B}(z)$}

\begin{definition}
For any  matrices $A$ and $B$ in $\Bbb{C}^{N\times N}$, $A+I$ and $B+I$ satisfy  (\ref{2.1}),
we define the Bateman matrix polynomials in the form
\begin{align}
&
\displaystyle
\mathbf{B}_{n}^{A,B}(z)=\ {}_{1}F_{2}(-nI;A+I,B+I;z)
\nonumber
\\
&=
\displaystyle
\sum_{k=0}^{\infty}\frac {z^{k}}{k!}(-nI)_{k}\left[(A+I)_{k}\right]^{-1}\left[(B+I)_{k}\right]^{-1}
\nonumber
\\
&=
\displaystyle
\sum_{k=0}^{\infty}\frac {(-1)^{k}n!z^{k}}{k!(n-k)!}
\left[(A+I)_{k}\right]^{-1}
\left[(B+I)_{k}\right]^{-1}.
\label{4.1}
\end{align}
\end{definition}

The initial conditions  are
\begin{align}
&
\displaystyle
\mathbf{B}_{0}^{A,B}(z)=I,
\nonumber
\\
&
\displaystyle
\mathbf{B}_{1}^{A,B}(z)=-(A+I)^{-1}(B+I)^{-1}z+I,
\nonumber
\\
&
\displaystyle
\mathbf{B}_{2}^{A,B}(z)=(A+I)^{-1}(A+2I)^{-1}(B+I)^{-1}(B+2I)^{-1}z^{2}-2(A+I)^{-1}(B+I)^{-1}z+I.
\nonumber
\end{align}

\begin{theorem}
\label{thm 4.1}
Let $A$, $B$ and $C$  be commutative matrices in $\Bbb{C}^{N\times N}$ satisfying (\ref{2.1}).
Then, the generating matrix functions for the Bateman matrix polynomials  are
\begin{eqnarray*}
\displaystyle
\sum_{n=0}^{\infty}\frac {t^{n}}{n!}\mathbf{B}_{n}^{A,B}(z)=e^{t}\ {}_{0}F_{2}(-;A+I,B+I;-zt),
\end{eqnarray*}
\begin{eqnarray*}
\displaystyle
\sum_{n=0}^{\infty}\frac {t^{n}}{n!}(C)_{n}\mathbf{B}_{n}^{A,B}(z)=(1-t)^{-C}\ {}_{1}F_{2}\left(C;A+I,B+I;-\frac {zt}{1-t}\right),
\
\displaystyle
|t|<1
\end{eqnarray*}
and
\begin{eqnarray*}
\displaystyle
\sum_{n=0}^{\infty}t^{n}\mathbf{B}_{n}^{A,B}(z)=e^{t}\ \mathbf{J}_{A,B}\left(3\sqrt[3]{zt}\right),
\end{eqnarray*}
where $\mathbf{J}_{A,B}\left(3\sqrt[3]{zt}\right)$ are hyper-Bessel matrix functions \cite{sa}.
\end{theorem}

\begin{theorem}
\label{thm 4.2}
Let $A$ and $B$ be commutative matrices in $\Bbb{C}^{N\times N}$ satisfying  (\ref{2.1}).
Then, the following recurrence relations hold
\begin{eqnarray*}
\displaystyle
\frac {d}{dz}\mathbf{B}_{n}^{A,B}(z)=n(A+I)^{-1}(B+I)^{-1}\mathbf{B}_{n-1}^{A+2I,B+2I}(z),
\
\displaystyle
n\geq1,
\end{eqnarray*}
\begin{eqnarray*}
\displaystyle
z\frac {d}{dz}\mathbf{B}_{n}^{A,B}(z)=n\mathbf{B}_{n}^{A,B}(z)-n\mathbf{B}_{n-1}^{A,B}(z),
\
\displaystyle
n\geq 1,
\end{eqnarray*}
\begin{eqnarray*}
\displaystyle
z\mathbf{B}_{n}^{A+I,B+I}(z)=\mathbf{B}_{n+1}^{A,B}(z)+(A+(n+1)I)\left(B+(n+1)I\right)
\mathbf{B}_{n}^{A,B}(z),
\end{eqnarray*}
\begin{eqnarray*}
\displaystyle
\mathbf{B}_{n+1}^{A,B}(z)=\mathbf{B}_{n+1}^{A+I,B}(z)+(n+1)\left(B+(n+1)I\right)\mathbf{B}_{n}^{A+I,B}(z)
\end{eqnarray*}
and
\begin{eqnarray*}
\displaystyle
\mathbf{B}_{n+1}^{A,B}(z)=\mathbf{B}_{n+1}^{A,B+I}(z)+(n+1)\left(A+(n+1)I\right)
\mathbf{B}_{n}^{A,B+I}(z).
\end{eqnarray*}
\end{theorem}

\begin{theorem}
\label{thm 4.3}
Let $A$ and $B$ be commutative matrices in $\Bbb{C}^{N\times N}$ satisfying  (\ref{2.1}).
Then, the following multiplication formulas hold true
\begin{eqnarray*}
\displaystyle
\mathbf{B}_{n}^{A,B}(zw)=\sum_{k=0}^{n}\binom{n}{k}w^{k}(1-w)^{n-k}\mathbf{B}_{k}^{A,B}(z)
\end{eqnarray*}
and
\begin{eqnarray*}
\displaystyle
\mathbf{B}_{n}^{A,B}(zw)=n!\left[(A+I)_{n}\right]^{-1}
\sum_{k=0}^{n}\left[(B+I)_{n}\right]^{-1}w^{k}L_{n-k}^{A+kI}(w)L_{k}^{B}(z).
\end{eqnarray*}
\end{theorem}

\begin{theorem}
\label{thm 4.4}
Let $A$ and $B$ be commutative matrices in $\Bbb{C}^{N\times N}$ satisfying  (\ref{2.1}).
Then, by virtue of the inversion formula, the following relation holds true
\begin{eqnarray*}
\displaystyle
z^{n}I=(A+I)_{n}(B+I)_{n}\sum_{k=0}^{n}\binom{n}{k}(-1)^{k}\mathbf{B}_{k}^{A,B}(z).
\end{eqnarray*}
\end{theorem}

\begin{theorem}
\label{thm 4.5}
Let $A+I$ and $B+I$ be commutative matrices in $\Bbb{C}^{N\times N}$ satisfying (\ref{2.1}).
Then, we obtain  the differential equation of third order
\begin{eqnarray}
\displaystyle
\left[z(\theta-n)I-\theta(\theta I+A)(\theta I+B)\right]\mathbf{B}_{n}^{A,B}(z)=\mathbf{0}.
\label{4.14}
\end{eqnarray}
\end{theorem}

\begin{theorem}
\label{thm 4.5}
Let $A$ and $B$ be commutative matrices in $\Bbb{C}^{N\times N}$ satisfying (\ref{2.1}).
Then, we obtain the third order differential equation
\begin{align}
\displaystyle
z^{2}\frac {d^{3}}{dz^{3}}\mathbf{B}_{n}^{A,B}(z)+(A+B+3I)z\frac {d^{2}}{dz^{2}}\mathbf{B}_{n}^{A,B}(z)
\left[(A+I)(B+I)-zI\right]
\frac {d}{dz}\mathbf{B}_{n}^{A,B}(z)+n\mathbf{B}_{n}^{A,B}(z)=\mathbf{0}.
\label{4.15}
\end{align}
\end{theorem}

\begin{proof}
Using the relations
\begin{align*}
&
\displaystyle
\theta(z^{n})=nz^{n},
\
\displaystyle
\theta^{2}=z^{2}\frac {d^{2}}{dz^{2}}+z\frac {d}{dz},
\\
&
\displaystyle
\theta^{3}=z^{3}\frac {d^{3}}{dz^{3}}+3z^{2}\frac {d^{2}}{dz^{2}}+z\frac {d}{dz}.
\end{align*}
and (\ref{4.14}), we obtain (\ref{4.15})
\end{proof}

\begin{theorem}
\label{thm 4.6}
Let $A$ and $B$ be commutative matrices in $\Bbb{C}^{N\times N}$ satisfying  (\ref{2.1}).
Then, we have
\begin{eqnarray*}
\displaystyle
\frac {\mathbf{B}_{n-k}^{A,B}(z)}{(n-k)!}=(-1)^{k}(\theta-n)_{k}\frac {\mathbf{B}_{n}^{A,B}(z)}{n!},
\
0\leq k \leq n,
\end{eqnarray*}
where the operator $ \theta $ is defined in the form
\begin{align}
&
\displaystyle
\theta=(\theta-n)+n=(\theta-n)_{1}+n,
\nonumber
\\
&
\displaystyle
\theta^{2}=(\theta-n)_{2}+(2n-1)(\theta-n)_{1}+n^{2},
\nonumber
\\
&
\displaystyle
\theta^{3}=(\theta-n)_{3}+3(n-1)(\theta-n)_{2}+\left(3n^{2}-3n+1\right)(\theta-n)_{1}+n^{3}.
\label{4.18}
\end{align}
\end{theorem}

\begin{theorem}
\label{thm 4.7}
Let $A$ and $B$ be commutative matrices in $\Bbb{C}^{N\times N}$ satisfying (\ref{2.1}).
Then, the following differential recurrence relation holds true
\begin{align}
&
\displaystyle
\bigg[n(A+nI)(B+nI)-\left(\left(3n^{2}-3n+1\right)I+(2n-1)(A+B)+AB-zI\right)(\theta-n)
\nonumber
\\
&
\qquad
\displaystyle
+\left(A+B+(3n-3)I\right)(\theta-n)_{2}I+(\theta-n)_{3}I\bigg]
\mathbf{B}_{n}^{A,B}(z)=\mathbf{0}.
\label{4.19}
\end{align}
\end{theorem}

\begin{proof}
(\ref{4.19}) follows easily from  (\ref{4.18}) and (\ref{4.14}).
\end{proof}

\begin{theorem}
\label{thm 4.7}
Let $A$ and $B$ be commutative matrices in $\Bbb{C}^{N\times N}$ satisfying  (\ref{2.1}).
Then, the following recurrence relation holds true
\begin{align}
&
\displaystyle
(A+nI)(B+nI)\mathbf{B}_{n}^{A,B}(z)
\nonumber
\\
&=
\displaystyle
\left(  3n^{2}I+(2A+2B-3I)n+(A-I)(B-I)-zI\right)
\mathbf{B}_{n-1}^{A,B}(z) - (n-1)\left(A+B+(3n-3)I\right)
\nonumber
\\
&
\quad
\displaystyle
\mathbf{B}_{n-2}^{A,B}(z)+(n-1)(n-2)\mathbf{B}_{n-3}^{A,B}(z),
\
\displaystyle
n\geq 3.
\label{4.20}
\end{align}
\end{theorem}

\begin{proof}
Put $\Psi_{n}^{A,B}(z)=\frac {\mathbf{B}_{n}^{A,B}(z)}{n!}$ in (\ref{4.1}), then
\begin{eqnarray*}
\displaystyle
\Psi_{n}^{A,B}(z)=\sum_{k=0}^{\infty}\frac {(-1)^{k}z^{k}}{k!(n-k)!}
\left[(A+I)_{k}\right]^{-1}
\left[(B+I)_{k}\right]^{-1}=\sum_{k=0}^{\infty}\phi(k,n;A,B;z).
\end{eqnarray*}
Now, we have
\begin{eqnarray*}
\displaystyle
\Psi_{n-1}^{A,B}(z)=\sum_{k=0}^{\infty}\frac {(-1)^{k}z^{k}}{k!(n-k-1)!}\left[(A+I)_{k}\right]^{-1}
\left[(B+I)_{k}\right]^{-1}=\sum_{k=0}^{\infty}(n-k)\phi(k,n;A,B;z),
\end{eqnarray*}
\begin{eqnarray*}
\displaystyle
\Psi_{n-2}^{A,B}(z)=\sum_{k=0}^{\infty}\frac {(-1)^{k}z^{k}}{k!(n-k-2)!}
\left[(A+I)_{k}\right]^{-1}
\left[(B+I)_{k}\right]^{-1}=\sum_{k=0}^{\infty}(n-k)(n-k-1)\phi(k,n;A,B;z)
\end{eqnarray*}
and
\begin{eqnarray*}
\displaystyle
\Psi_{n-3}^{A,B}(z)=\sum_{k=0}^{\infty}\frac {(-1)^{k}z^{k}}{k!(n-k-3)!}
\left[(A+I)_{k}\right]^{-1}
\left[(B+I)_{k}\right]^{-1}=\sum_{k=0}^{\infty}(n-k)(n-k-1)(n-k-2)\phi(k,n;A,B;z).
\end{eqnarray*}
Also
\begin{align*}
\displaystyle
\Psi_{n-1}^{A,B}(z)
&=
\displaystyle
\sum_{k=0}^{\infty}\frac {(-1)^{k}z^{k+1}}{k!(n-k-1)!}
\left[(A+I)_{k}\right]^{-1}
\left[(B+I)_{k}\right]^{-1}
\nonumber
\\
&=
\displaystyle
\sum_{k=1}^{\infty}\frac {(-1)^{k-1}z^{k}}{(k-1)!(n-k)!}
\left[(A+I)_{k-1}\right]^{-1}
\left[(B+I)_{k-1}\right]^{-1}
\nonumber
\\
&=
\displaystyle
\sum_{k=0}^{\infty}(-k)(A+kI)(B+kI)\phi(k,n;A,B;z).
\end{align*}
Then there exists the relation
\begin{eqnarray*}
\displaystyle
\Psi_{n}^{A,B}(z)+(a+bz)\Psi_{n-1}^{A,B}(z)+c \Psi_{n-2}^{A,B}(z)+d \Psi_{n-3}^{A,B}(z)=0,
\end{eqnarray*}
where $a$, $b$, $c$ and $d$ are constants dependent only on $n$.

We are led to the identity in $k$
\begin{eqnarray*}
\displaystyle
I+a(n-k)I-bk(A+kI)(B+kI)+c(n-k)(n-k-1)I+d(n-k)(n-k-1)(n-k-2)I=0.
\end{eqnarray*}
By equating coefficients, we find  $b=\frac {1}{n}(A+nI)^{-1}(B+nI)^{-1}$,
$d=-\frac {1}{n}(A+nI)^{-1}(B+nI)^{-1}$,
$a=\frac {AB+(2n-1)(A+B)+\left(3n^{2}-3n+1\right)I}{n}(A+nI)^{-1}(B+nI)^{-1}$
and  $c=\frac {(n-1)\left(A+B+(3n-3)I\right)}{n(n-1)}(A+nI)^{-1}(B+nI)^{-1}$.
We obtain
\begin{align*}
&
\displaystyle
n(A+nI)(B+nI)\frac {\mathbf{B}_{n}^{A,B}(z)}{n!}-\left[AB+(2n-1)(A+B)+(3n^{2}-3n+1)I-xI\right]
\frac {\mathbf{B}_{n-1}^{A,B}(z)}{(n-1)!}
\nonumber
\\
&
\quad
\displaystyle
+\left(A+B+(3n-3)I\right)\frac {\mathbf{B}_{n-2}^{A,B}(z)}{(n-2)!}-\frac {\mathbf{B}_{n-3}^{A,B}(z)}{(n-3)!}=0
\end{align*}
or
\begin{align*}
&
\displaystyle
(A+nI)(B+nI)\mathbf{B}_{n}^{A,B}(z)-\left[AB+(2n-1)(A+B)+\left(3n^{2}-3n+1\right)I-xI\right]\mathbf{B}_{n-1}^{A,B}(z)
\nonumber
\\
&
\quad
\displaystyle
+(n-1)(A+B+(3n-3)I)\mathbf{B}_{n-2}^{A,B}(z)-(n-1)(n-2)\mathbf{B}_{n-3}^{A,B}(z)=0.
\end{align*}
Hence, the theorem.
\end{proof}

\subsection{Bateman matrix polynomials $\mathbf{J}_{n}^{A,B}(z)$}

\begin{definition}
Let $A$ and $B$ be matrices in $\Bbb{C}^{N\times N}$, $A+I$ and $\frac {1}{2}A+B+I$ satisfy (\ref{2.1}).
We define the Bateman matrix polynomials $\mathbf{J}_{n}^{A,B}(z)$ by
\begin{eqnarray*}
\displaystyle
\mathbf{J}_{n}^{A,B}(z)=\frac {\Gamma\left(\frac {1}{2}A+B+nI+I\right)}{n!}
\Gamma^{-1}(A+I)\Gamma^{-1}\left(\frac {1}{2}A+B+I\right)z^{A}\ {}_{1}F_{2}\left(-nI;A+I,B+\frac {1}{2}A+I;z^{2}\right).
\end{eqnarray*}
\end{definition}

\begin{corollary}
\label{thm 4.8}
For the   Bateman matrix polynomials, we have
\begin{eqnarray}
\displaystyle
\mathbf{J}_{n}^{A,B-\frac {1}{2}A}\left(\sqrt{z}\right) =
\frac {\Gamma(B+nI+I)}{n!}\Gamma^{-1}(A+I)\Gamma^{-1}(B+I)z^{\frac {1}{2}A}\mathbf{B}_{n}^{A,B}(z).
\label{4.33}
\end{eqnarray}
\end{corollary}

\begin{theorem}
\label{thm 4.9}
The Bateman matrix polynomials $\mathbf{J}_{n}^{A,B}(z)$ have the following recurrence relations
\begin{eqnarray}
\displaystyle
J_{n}^{A,B}(z)=J_{n-1}^{A,B}(z)+J_{n}^{A,B-I}(z),
\label{4.34}
\end{eqnarray}
\begin{eqnarray}
\displaystyle
AJ_{n}^{A,B-\frac {1}{2}I}(z)=zJ_{n}^{A-I,B}(x)+zJ_{n-1}^{A+I,B}(z)
\label{4.35}
\end{eqnarray}
and
\begin{eqnarray}
\displaystyle
(A+nI)J_{n}^{A,B}(z)=\left(nI+B+\frac {1}{2}A\right)J_{n-1}^{A,B}(z)+zJ_{n}^{A-I,B+\frac {1}{2}I}(z).
\label{4.36}
\end{eqnarray}
\end{theorem}

\begin{proof}
An an illustration of the proof,  set $A=-nI$, $B=A+I$, $C=\frac {1}{2}A+B+I$ and $z=z^{2}$ in (\ref{2.20}),
we obtain the contiguous function relations (\ref{4.34})-(\ref{4.36}).
\end{proof}

\begin{theorem}
\label{thm 4.10}
The functions $J_{n}^{A,B}(x)$ and $\mathbf{B}_{n}^{A,B}(z)$ satisfy the following  recurrence relation
\begin{align}
&
\displaystyle
n(B+nI)J_{n}^{A,B}(z)
\nonumber
\\
&=
\displaystyle
\left[\left(B+nI+\frac {1}{2}A\right)(A+nI)+(n-1)\left(B+\frac {3}{2}A+2nI-I\right)-z^{2}I\right]J_{n-1}^{A,B}(z)
\nonumber
\\
&
\quad
\displaystyle
\left(B+\frac {1}{2}A+nI-I\right)
\left(B+\frac {3}{2}A+3n-3I\right)J_{n-2}^{A,B}(z)
\nonumber
\\
&
\quad
\displaystyle
+\left(B+\frac {1}{2}A+nI-I\right)
\left(B+\frac {1}{2}A+nI-2I\right)J_{n-3}^{A,B}(z),
\
\displaystyle
n\geq3.
\label{4.37}
\end{align}
\end{theorem}

\begin{proof}
Using (\ref{4.20}) and (\ref{4.33}), we obtain (\ref{4.37}).
\end{proof}

\subsection{Young's matrix functions}

\begin{definition}
Let $A$ be a  matrix in $\Bbb{C}^{N\times N}$, $\frac {1}{2}(A+I)$ and $\frac {1}{2}(A+2I)$ satisfying (\ref{2.1}).
We define the Young's matrix function by
\begin{eqnarray*}
\displaystyle
Y_{A}(x)=z^{A}\sum_{k=0}^{\infty}(-1)^{k}\Gamma^{-1}\left(A+(2k+1)I\right)x^{2k}=
x^{A}\Gamma^{-1}(A+I)\ {}_{1}F_{2}\left(I;\frac {1}{2}(A+I),\frac {1}{2}(A+2I);-\frac {x^{2}}{4}\right).
\end{eqnarray*}
\end{definition}

\begin{theorem}
The following relations for the Young's matrix functions hold true
\begin{eqnarray}
\displaystyle
Y_{A}(x)=\sum_{k=0}^{\infty}\frac {1}{k!}(2z)^{\frac {1}{2}(A+2kI)}
x^{2k}\left(\frac {1}{2}(A-I)\right)_{k}
\left[(A+I)_{2k}\right]^{-1}\Gamma(A+I)^{-1}
\Gamma\left(\frac {1}{2}(A+(2k+2)I)\right)
J_{\frac {1}{2}(A+2kI)}(x)
\label{4.39}
\end{eqnarray}
and
\begin{eqnarray}
\displaystyle
Y_{A}(x)=\frac {1}{2}\sum_{k=0}^{\infty}\frac {1}{k!}
\left(\frac {1}{2}A\right)_{k}
\left[(A+I)_{2k}\right]^{-1}(2x)^{\frac {1}{2}\left(A+(2k+1)I\right)}
\Gamma^{-1}(A+I)\Gamma\left(\frac {1}{2}\left(A+(2k+1\right)I)\right)J_{\frac {1}{2}\left(A+(2k-1)I\right)}(x),
\label{4.40}
\end{eqnarray}
where $J_{A}(z)$ are the Bessel matrix functions   \cite{sj}
\begin{eqnarray*}
\displaystyle
J_{A}(x)=\Gamma^{-1}(A+I)\left(\frac {x}{2}\right)^{A}\ {}_{0}F_{1}\left(-;A+I;-\frac {x^{2}}{4}\right),
\end{eqnarray*}
\begin{eqnarray}
\displaystyle
J_{\frac {1}{2}A}(x)=\sum_{k=0}^{\infty}\frac {(-1)^{k}}{k!}
(2x)^{-\frac {1}{2}A}\left(\frac {1}{2}(A-I)\right)_{k}[(A+I)_{2k}]^{-1}
\Gamma^{-1}\left(\frac {1}{2}(A+2I)\right)
\Gamma\left(A+(2k+1)I\right)Y_{A+2kI}(x)
\label{4.41}
\end{eqnarray}
and
\begin{eqnarray}
\displaystyle
J_{\frac {1}{2}(A-I)}(x)=2\sum_{k=0}^{\infty}\frac {(-1)^{k}}{k!}\left(\frac {1}{2}A\right)_{k}
(2x)^{-\frac {1}{2}A}\left[(A+I)_{2k}\right]^{-1}\Gamma^{-1}\left(\frac {1}{2}(A+I)\right)\Gamma\left(A+(2k+1)I\right)Y_{A+2kI}(x).
\label{4.42}
\end{eqnarray}
\end{theorem}

\begin{proof}
Putting $A=I$, $B=\frac {1}{2}(A+I)$, $C=\frac {1}{2}(A+2I)$ and $z=-\frac {x^{2}}{4}$ in (\ref{2.28})-(\ref{2.32})
and  multiplying  by  $z^{A}$, we obtain (\ref{4.39})-(\ref{4.42}).
\end{proof}

\begin{theorem}
The following  matrix differential equation for the Young's matrix function holds true
\begin{eqnarray*}
\displaystyle
x\frac {d^{3}Y_{A}(x)}{dx^{3}}+\left(2I-A\right) \frac {d^{2}Y_{A}(x)}{dx^{2}}+x\frac {dY_{A}(x)}{dx}+
\left(2I-A\right) Y_{A}(x)=\mathbf{0}.
\end{eqnarray*}
\end{theorem}

\begin{proof}
In (\ref{2.23}), set $A=I$, $B=\frac {1}{2}(A+I)$  and $C=\frac {1}{2}(A+2I)$.
We can rewrite the matrix differential equation
\begin{eqnarray*}
\displaystyle
z^{2}\frac {d^{3}}{dz^{3}}\ {}_{1}F_{2}+\frac {1}{2}(2A+5I)z\frac {d^{2}}{dz^{2}}\ {}_{1}F_{2}+
\left(\frac {1}{4}(A+I)(A+2I)-z I\right)\frac {d}{dz}\ {}_{1}F_{2}-\ {}_{1}F_{2}=\mathbf{0}.
\end{eqnarray*}
As an illustration, set $z=-\frac {x^{2}}{4}$ in (\ref{2.23}), therefore
\begin{eqnarray*}
\displaystyle
\frac {dW}{dz}=-\frac {2}{x}\frac {dW}{dx},
\end{eqnarray*}
\begin{eqnarray*}
\displaystyle
\frac {d^{2}W}{dz^{2}}=\frac {4}{x^{2}}\frac {d^{2}W}{dx^{2}}-\frac {4}{x^{3}}\frac {dW}{dx},
\end{eqnarray*}
\begin{eqnarray*}
\displaystyle
\frac {d^{3}W}{dz^{3}}=-\frac {8}{x^{3}}\frac {d^{3}W}{dx^{3}}+\frac {24}{x^{4}}\frac {d^{2}W}{dx^{2}}-\frac {24}{x^{5}}\frac {dW}{dx},
\end{eqnarray*}
and we obtain
\begin{eqnarray}
\displaystyle
x^{2}\frac {d^{3}W}{dx^{3}}+2\left(A+I\right)x\frac {d^{2}W}{dx^{2}}+\left(x^{2}I+A(A+I)\right)\frac {dW}{dx}+2xW=\mathbf{0}.
\label{4.45}
\end{eqnarray}
One solution (\ref{4.45}) is
\begin{eqnarray*}
\displaystyle
W=\ {}_{1}F_{2}\left(I;\frac {1}{2}(A+I),\frac {1}{2}(A+2I);-\frac {x^{2}}{4}\right).
\end{eqnarray*}
Note that
\begin{eqnarray*}
\displaystyle
W(x)=\Gamma(A+I)x^{-A}Y_{A}(x),
\end{eqnarray*}
\begin{eqnarray*}
\displaystyle
\frac {dW(x)}{dx}=\Gamma(A+I)\left[x^{-A}\frac {dY_{A}(x)}{dx}-Ax^{-A-I}Y_{A}(x)\right],
\end{eqnarray*}
\begin{eqnarray*}
\displaystyle
\frac {d^{2}W(x)}{dx^{2}}=\Gamma(A+I)\left[x^{-A}\frac {d^{2}Y_{A}(x)}{dx^{2}}-2Ax^{-A-I}\frac {dY_{A}(x)}{dx}+A(A+I)x^{-A-2I}Y_{A}(x)\right],
\end{eqnarray*}
and
\begin{align*}
\displaystyle
\frac {d^{3}W(x)}{dx^{3}}
&=
\displaystyle
\Gamma(A+I)\bigg[x^{-A}\frac {d^{3}Y_{A}(x)}{dx^{3}}-3Ax^{-A-I}\frac {d^{2}Y_{A}(x)}{dx^{2}}
\\
&
\quad
\displaystyle
3A(A+I)x^{-A-2I}\frac {dY_{A}(x)}{dx}-A(A+I)(A+2I)x^{-A-3I}Y_{A}(x)\bigg].
\end{align*}
We obtain
\begin{align}
&
\displaystyle
x^{3}\frac {d^{3}Y_{A}(x)}{dx^{3}}-3Ax^{2}\frac {d^{2}Y_{A}(x)}{dx^{2}}+3A(A+I)x\frac {dY_{A}(x)}{dx}-A(A+I)(A+2I)Y_{A}(x)
\nonumber
\\
&
\quad
\displaystyle
+2(A+I)\left[x^{2}\frac {d^{2}Y_{A}(x)}{dx^{2}}-2Ax\frac {dY_{A}(x)}{dx}+A(A+I)Y_{A}(x)\right]
\nonumber
\\
&
\quad
\displaystyle
+\left(x^{2}I+A(A+I)\right)\left[x\frac {dY_{A}(x)}{dx}-AY_{A}(x)\right]+2x^{2}Y_{A}(x)=\mathbf{0}.
\label{4.46}
\end{align}
We seek an equation satisfied by $W(x)=\Gamma(A+I)x^{-A}Y_{A}(x)$.
In (\ref{4.46}), we  put $W(x)=\Gamma(A+I)x^{-A}Y_{A}(x)$ and arrive at the matrix differential equation
\begin{eqnarray}
\displaystyle
x^{3}\frac {d^{3}Y_{A}(x)}{dx^{3}}+(2I-A)x^{2}\frac {d^{2}Y_{A}(x)}{dx^{2}}+x^{3}\frac {dY_{A}(x)}{dx}+
(2I-A)x^{2}Y_{A}(x)=\mathbf{0},
\label{4.47}
\end{eqnarray}
of which one solution is $Y_{A}(x)$ in (\ref{4.47}).
\end{proof}

\section{Conclusions}

This paper has delved into an exploration of Bateman's matrix polynomials with integral indices.
Throughout our investigation, we have uncovered a range of significant findings.
These include the identification of key properties related to hypergeometric matrix functions,
the establishment of integral representations for generalized hypergeometric matrix functions,
and the derivation of crucial recurrence relations.
Furthermore, our study has introduced novel concepts such as matrix differential equations,
integral transforms, and fractional integral formulas for generalized gauss hypergeometric functions.
These innovations have been achieved through the adept utilization of various mathematical tools
such as the beta and Laplace transform formulas, Mellin transforms of integrable functions,
Erdelyi Kober-type fractional integral operators, as well as fractional integrals and derivatives
pertaining to  Riemann–Liouville and Weyl operators.
Ultimately, these theoretical advancements can  find practical applications,
particularly in the realm of Bateman's matrix polynomials across diverse types and Young's matrix functions.
This research opens avenues for further exploration and development within this intricate field of study.

\subsection*{Ethics approval}

Not applicable.

\subsection*{Funding}

Not applicable.

\subsection*{Conflict of interest}

All of the   authors have no conflicts of interest.

\subsection*{Data availability statement}

Not applicable.

\subsection*{Code availability}

Not applicable.

\subsection*{Consent to participate}

Not applicable.

\subsection*{Consent for publication}

Not applicable.

\end{document}